\documentclass[a4paper,12pt,oneside]{article}
\usepackage[english]{babel}
\usepackage[T1]{fontenc} 
\usepackage[utf8]{inputenc}
\usepackage{amsthm}
\usepackage{amsmath}
\usepackage{amssymb}  
\usepackage{indentfirst}
\usepackage{fancyhdr}
\usepackage{amsthm}
\usepackage{graphicx}
\usepackage{pdfpages}
\usepackage{url}
\usepackage{hyperref}

\definecolor{dg}{rgb}{0.01, 0.75, 0.24}

\numberwithin{equation}{section}

\theoremstyle{plain} 
\newtheorem{thm}{Theorem}[section]

\newtheorem{lem}[thm]{Lemma} 
\newtheorem{prop}[thm]{Proposition} 
\newtheorem{rmk}[thm]{Remark}

\def\it2{\int_{B_{\tau_2}}}

\def\ir {\int_{B_{r}}}

\def\ibr2 {\int_{B_{\frac{R}{2}}}}
\def\ibro2 {\int_{B_{\frac{\rho}{2}}}}

\def\fhi{\varphi}
\def\R{\mathbb{R}}
\def\N{\mathbb{N}}

\def\dd{\textrm{d}}

\def\fibr2 {\displaystyle\fint_{B_{\frac{R}{2}}}}
\def\fibro2 {\fint_{B_{\frac{\rho}{2}}}}

\usepackage{geometry}
\allowdisplaybreaks[4]

\begin{document}

\title{\textbf{Regularity results for H\"{o}lder minimizers
to functionals with non-standard growth
}}

\author {   
\sc{Antonio Giuseppe Grimaldi}\thanks{Dipartimento di Scienze Matematiche "G. L. Lagrange",
Politecnico di Torino, Corso Duca degli Abruzzi, 24, 10129 Torino, 
 Italy. E-mail: \textit{antonio.grimaldi@polito.it}} 
  \;\textrm{and}   Erica Ipocoana \thanks{Freie Universit\"at Berlin,
Department of Mathematics and Computer Science, Arnimallee 9, 
14195 Berlin, Germany.
E-mail: \textit{erica.ipocoana@fu-berlin.de}}}

\maketitle
\maketitle

\begin{abstract}
We study the regularity properties of H\"older continuous minimizers to non-autonomous functionals satisfying $(p,q)$-growth conditions, under Besov assumptions on the coefficients. In particular, we are able to prove higher integrability and higher differentiability results for solutions to our minimum problem.
\end{abstract}


\section{Introduction}
In this work, we aim to present some regularity results for the solutions to the problem
\begin{equation}\label{minprob}
\min\left\lbrace \mathcal{F}(u,\Omega):=\int_\Omega F(x,Du) \  \mathrm{d}x: u\in W^{1,p}(\Omega, \R^N) 
\right\rbrace, 
\end{equation}
where $N\geq 1$, $\Omega$ is a bounded open subset of $\R^n$, for $n\geq 2$, and
$F:\Omega\times\R^{N\times n}\to [0,+\infty)$ is a Carath\'{e}odory function. 
Besides, we assume that there exist positive constants $\nu, L$ and exponents $1< p < q < +\infty$ such that

$$ \nu|z|^p \leq F(x,z) \leq L(1+|z|^q)\eqno{\rm{{ (F1)}}} $$

$$\nu(\mu^2+|z_1|^2+|z_2|^2)^\frac{p-2}{2} \leq \langle \partial_z F(x, z_1)-\partial_z F(x, z_2),z_1-z_2\rangle   \eqno{\rm{{ (F2)}}}$$

\noindent for a.e. $x \in \Omega$ and all $z,z_1,z_2 \in \mathbb{R}^{N \times n}$, with $\mu \in [0,1]$.
Note that the ellipticity assumption (F2) and (F1) imply 
$$ |\partial_z F(x,z)| \leq c  (1 +|z|^{q-1}) \eqno{\rm{{ (F3)}}}$$

\noindent for a.e. $x \in \Omega$ and every $z \in \mathbb{R}^{N \times n}$, with a constant $c:=c(L,q)$. {Moreover, we assume that
there exists a non-negative function 
$g \in L^{\sigma}_{\text{loc}}(\Omega)$,
where $\sigma$ is finite and its precise expression will be introduced in \eqref{sigma}, such that}
{$$ |\partial_z F(x,z)-\partial_z F(y, z)| \leq |x-y|^{\alpha} (g(x)+g(y)) (1 +|z|^{q-1}) \eqno{\rm{{ (F4)}}}$$
\noindent for $\alpha \in (0,1)$, a.e. $x,y \in \Omega$ and every $z \in \mathbb{R}^{N \times n}$.\\ Note that, according to standard results recalled in Section \ref{secbesov} and in particular Theorem \ref{thm2.7}, assumption (F4) states that the function $x \mapsto \partial_z F(x,z)$ belongs to $B_{\sigma,\infty}^{\alpha}$.
}

Functionals satisfying assumption  $\textrm{(F1)}$ above are called $(p,q)$-growth functionals, as first referred by Marcellini \cite{marcellini1986, marcellini1989, marcellini1991}. As a necessary and sufficient condition for regularity of minima, these non-standard growth functionals satisfy gap bounds of the type
$$\dfrac{q}{p}<1+o(n), \qquad o(n)\approx \dfrac{1}{n},$$
according to \cite{25dfm,marcellini1991}. 

{A very interesting model case of functionals with $(p,q)$-growth is the so called double phase, whose energy density has the form 
\begin{equation}\label{dphase}
    F(x,z):=|z|^p+a(x)|z|^q
\end{equation}
with $a(x) \in \mathcal{C}^{0,\alpha}(\Omega)$. This model case, first introduced by Zhikov \cite{60dfm,61dfm,62dfm} in order to describe the behaviour of strongly anisotropic materials in the context of homogenisation, nonlinear elasticity and Lavrentiev phenomenon, has been widely investigated in \cite{BCM,25dfm,FMM} with the aim of identifying sufficient and necessary conditions on the relation between $p$, $q$ and $\alpha$ to establish the H\"older continuity of the gradient of the local minimizers. The sharp condition on the gap is given by
$$\dfrac{q}{p} \le 1+ \dfrac{\alpha}{n}. $$
On the other hand, in \cite{CM} the authors showed that assuming that the minimizers are a priori bounded yields a gap bound of the form
\begin{equation*}
    q \le p + \alpha,
\end{equation*}
that is sharp as shown in \cite{25dfm,FMM}.
This is independent of the dimension and this phenomenon holds true for a priori bounded minimizers for functionals with $(p,q)$-growth without the structure as \eqref{dphase} (see for instance \cite{BCM,carozza0,choe, grimaldi}). 
\\The fact that the a priori boundedness of minima leads to non-dimensional bounds on
the distance $q-p$ has been the motivation to the study of a priori more regular minimizers. 
}
{Actually, in \cite{CM} it has been proved that in the case of the double phase functional
the gap improves with the a priori $\mathcal{C}^{0,\gamma}$-regularity of the minima, for $0<\gamma<1$. Indeed, the gap bound takes the form
\begin{equation*}
q<p+\dfrac{\alpha}{1-\gamma},
\end{equation*}
that is still independent of the dimension $n$ and is sharp, as shown in \cite{2dfm}.}
Moreover, for $\gamma\to 1$, it holds
\begin{equation}\label{boundO}
q<p+o(\gamma), \qquad o(\gamma)\approx \dfrac{1}{1-\gamma},
\end{equation}
which is in particular coherent with the fact that Lipschitz continuity of minima is a main issue for $(p,q)$-growth problems. {An additional} step regarding a priori H\"{o}lder continuous minimizers is made in \cite{DFM}, where the authors are able to improve the regularity of minima of non-autonomous functionals from $W^{1,p}$ to $W^{1,q}$, taking into account the bound \eqref{boundO}. {Here, we further extend the results obtained in \cite{DFM} by considering weaker assumptions on the function $g$, in particular, when $g$ is not bounded (see (F3)).}
\\

In the recent paper \cite{grimaldi}, it has been proved that a priori bounded solutions to obstacle problems with non-standard growth satisfy extra fractional differentiability properties, under a free dimensional condition between the ellipticity and the growth exponents of the functional. Here, we continue the study of regularity theory of a priori more regular minimizers to integral functionals of Calculus of Variations with $(p,q)$-growth.
The novelty of this work consists in showing that, assuming a Besov regularity on the coefficients of the problem, the higher integrability and differentiability of a priori H\"older continuous minima hold under a relation between $p$ and $q$ of type \eqref{boundO}. 
\\
For precise definition and properties of Besov spaces we refer to Section \ref{secbesov}.
We note that, indeed, recently, there have been many research activities regarding the analysis of higher differentiability (of integer and fractional order) of solutions to variational problems; see \cite{ambrosio1,ambrosio2,baison.clop2017,clop,Giova,Giova.Pass,Pass1,Pass2,Pass3}, for the unconstrained case, and \cite{byun1,caselli,defilippis1,defilippis,eleuteri.passarelli,eleuteri.passarelli1,gavioli1,gavioli2,Gentile,gentile0,gentile1,grimaldi0,grimaldi,grimaldi.ipocoana,grimaldi.ipocoana1,zhang.zheng}, for the obstacle problem.\\

As stated above, according to \eqref{minprob}, we are here dealing with a non-autonomous functional. Therefore, the Lavrentiev phenomenon has to be taken into account, 
i.e. 
$$\inf_{w \in W^{1,p}(\Omega)} \int_{\Omega} F(x,Dw) \ \mathrm{d} x 
< \inf_{w \in W^{1,q}(\Omega)} \int_{\Omega} F(x,Dw) \ \mathrm{d}x .$$
This is a clear obstruction to regularity, since it prevents minimizers to belong to $W^{1,q}$ (see \cite{marcellini1986, 60dfm, 62dfm} for more details).\\
Here, we are interested in the case of a priori H\"{o}lder minimizers to the problem \eqref{minprob}. Hence, we will adopt suitable definitions of relaxed and gap functionals, aiming at exploiting this fact. Let us fix constants $H>0$ and $\gamma \in (0,1)$.
Let $B\subset\Omega$ be a ball, we can consider the relaxed functional
\begin{equation}
\mathcal{\bar{F}}_{H,\gamma}(u,B) := \inf_{\lbrace u_j\rbrace\in \mathcal{C}_{H,\gamma}(u,B)}
\left\lbrace \liminf_{j\to\infty}\int_B F(x,Du_j)\ \mathrm{d}x\right\rbrace \label{lavrentiev}
\end{equation}
for any $u\in W^{1,1}(B,\R^N)$, where 
\begin{align}
\mathcal{C}_{H,\gamma}(u,B):=
\left\lbrace \lbrace w_j\rbrace \subset W^{1,\infty}(B,\R^N):w_j\rightharpoonup u \textrm{ in } W^{1,p}(B,\R^N),\ \sup_j [w_j]_{0,\gamma;B}\leq H \right\rbrace
\label{CHgamma}
\end{align}
and where we used the standard notation
 $$[v]_{0,\gamma;A} := \sup_{x,y \in A, x \neq y} \frac{|v(x)-v(y)|}{|x-y|^{\gamma}} ,$$
whenever $A \subset \R^n$ and $v : A \rightarrow \R^N$.  By lower semicontinuity of the functional $\mathcal{F}$ with respect to the weak convergence in $W^{1,p}$, we have $\tilde{\mathcal{F}}_{H,\gamma}(u,B) \ge \mathcal{F}(u,B)$. 
 We define for all $u \in W^{1,1}(B,\R^N)$ the Lavrentiev gap functional relative to $\mathcal{F}$ as follows
\begin{equation}\label{gapfun}
\mathcal{L}_{F,H,\gamma}(u,B):=
\begin{cases}
\bar{\mathcal{F}}_{H,\gamma}(u,B)-\mathcal{F}(u,B) \quad & \text{if} \  \mathcal{F}(u,B)< + \infty \\
 0 \quad & \text{if} \  \mathcal{F}(u,B)= + \infty.
\end{cases}
\end{equation}

From the previous definitions, we can easily derive 
\begin{prop}
    Let $w \in W^{1,1}(B,\R^N)$ be such that $\mathcal{F}(w,B)$ is finite, for some ball $B \Subset \Omega$. Then, $\bar{\mathcal{F}}_{H,\gamma}(w,B)=\mathcal{F}(w,B)$, for some $H>0$ if, and only if, there exists a sequence
$\lbrace w_j\rbrace\in \mathcal{C}_{H,\gamma}(w,B)$ such that
$$ \mathcal{F}(w_j,B) \rightarrow \mathcal{F}(w,B).$$
\end{prop}

In particular, in the sequel we need to impose that the Lavrentiev phenomenon does not occur, as a necessary condition for the the regularity of minima \cite{1dfm,20dfm,25dfm}.

We may now state the two main results contained in this work, concerning higher integrability and higher differentiability of solutions to \eqref{minprob}, respectively.
\begin{thm}[Higher integrability]\label{thmhi}
Let $\alpha, \gamma \in (0,1)$ and let $F(x,z)$ satisfy (F1) and (F2) for exponents $1<p<q$ such that 

\begin{equation}\label{gapqpx}
q<p+\dfrac{\min \lbrace \alpha,2\gamma\rbrace}{\Theta(1-\gamma)},
\end{equation} 
where 
\begin{eqnarray}\label{teta}
 \Theta:=
  \begin{cases}
  1 & \quad \text{if} \  p\geq 2  \\
  \dfrac{2}{p} & \quad \text{if} \  1<p<2.
  \end{cases} 
\end{eqnarray}
Moreover, assume that there exists a non-negative function $g \in L^{\sigma}_{\text{loc}}(\Omega)$, with
\begin{equation}\label{sigma}
\sigma = \dfrac{p(1-\gamma)+\xi}{(p-q)(1-\gamma)+\xi}
\end{equation}
and $(q-p)(1-\gamma)<\xi< {\dfrac{\min \lbrace \alpha,2\gamma\rbrace}{\Theta}}$,
{such that (F4) holds.}\\

Let $u\in W^{1,p}(\Omega, \mathbb{R}^N)$ be a solution to the minimum problem \eqref{minprob}.
Assuming that
\begin{equation}
    \mathcal{L}_{F,H,\gamma}(u,B_r)=0,\label{lav}
\end{equation}
holds for a ball $B_r \Subset \Omega$ with $r \le 1$ and for some $H >0$,
then
 \begin{equation}\label{hi}
 Du\in L^{\tilde{q}}_{\textrm{loc}}(\Omega, \R^{N \times n}),
 \end{equation} 
 for all $q<\tilde{q}<p+\dfrac{\min \lbrace \alpha,2\gamma\rbrace}{\Theta(1-\gamma)}$.
In particular, for every $B_\rho \Subset B_r$, we have
\begin{equation}
\Vert Du\Vert_{L^{\tilde{q}}(B_\rho)}\leq \dfrac{C}{(r-\rho)^{\kappa_1}}([\mathcal{F}(u,B_r)]^\frac{1}{p}+\Vert g \Vert^\sigma_{L^\sigma(B_r)}+H+1)^{\kappa_2},
\end{equation}
for constants $C:=C(n,p,q,\tilde{q}, \nu ,L, \alpha,\gamma,\sigma )$ and $\kappa_1,\kappa_2:=\kappa_1,\kappa_2 (n,p,q,\tilde{q},  \alpha,\gamma,\sigma )$.
\end{thm}

\begin{rmk} 
  Assumption \eqref{lav} in Theorem \ref{thmhi} ensures that the minimizer $u$ is H\"older continuous. By minimality
  we have $\mathcal{F}(u,B) < + \infty$, for every ball $B \Subset \Omega$. Note that $\mathcal{C}_{H,\gamma}(w,B)$ cannot be empty, since it would imply 
  $\bar{\mathcal{F}}(u,B) =+ \infty$. Therefore,
  $\mathcal{C}_{H,\gamma}(w,B) \neq \varnothing$, which yields $u \in \mathcal{C}^{0,\gamma}_{loc}(\Omega,\R^N)$, with $ [u]_{0,\gamma;B} \le H $. 
  \\ Besides, we point out that the gap condition \eqref{gapqpx} is of type \eqref{boundO}.
\end{rmk}

{It is convenient to introduce the auxiliary function
\begin{center}
$V_{p}(z)=(\mu^{2}+|z|^{2})^\frac{p-2}{4} z$
\end{center}
defined for all $z\in \mathbb{R}^{k}$, $k \in \N$, in order to state the following result and as it will be largely used in the sequel.}

\begin{thm}[Higher differentiability]\label{thmhd}
Let $\alpha,\gamma \in (0,1)$ and let $\sigma$ be the exponent defined in \eqref{sigma}. Let $F(x,z)$ satisfy (F1)--(F4) for exponents $1<p<q$ such that \eqref{gapqpx} holds.
Let $u\in W^{1,p}(\Omega,\R^N)$ be a solution to the minimum problem \eqref{minprob}.
Assuming that
\eqref{lav} holds for a ball $B_r \Subset \Omega$ with $r \le 1$ and for some $H >0$,
then
\begin{equation}
V_p(Du)\in B^\frac{\alpha}{2}_{2,\infty,\textrm{loc}}(\Omega, \R^{N \times n})
\end{equation}
and the following estimate
\begin{align}
 \displaystyle\int_{B_{\rho}}  |\tau_{h}V_p(Du)|^{2} \ \mathrm{d} x 
    \le C |h|^\alpha \biggl\{
    \displaystyle\int_{B_{r}}
     (1+|Du|^p) \ \mathrm{d}x +
    \Vert g \Vert^\sigma_{L^\sigma(B_{r})} + H
    \biggr\}^\vartheta
\end{align}
holds true for every $B_\rho \Subset B_r$, for constants $C:=C(n,p,q, \nu ,L, \alpha,\gamma,\sigma,r , \rho )$ and
$\vartheta:=\vartheta (n,p,q,  \alpha,\gamma,\sigma )$.
\end{thm}


The strategy of the proof of Theorem \ref{thmhi} is to approximate the original minimizer $u$ of \eqref{minprob} with a sequence $\{ u_j \}$ of $W^{1,q}$-regular solutions to suitable approximating problems. The possibility of this approximation relies on the assumptions \eqref{lav}. The main difficulty of this work is to build the approximating minimizers $\{ u_j \}$ such that they are uniformly bounded in $\mathcal{C}^{0,\gamma}$. Actually, we build the sequence $\{ u_j \}$ in a way that it is bounded in a fractional Sobolev space $W^{s,2d}$, which is a good approximation of $\mathcal{C}^{0,\gamma}$, for a suitable choice of parameters $s$ and $d$. In order to deal with this issue, we follow the strategy first proposed by De Filippis and Mingione in \cite{DFM}, i.e. we employ a method based on approximating the original local functional with mixed local and non-local functionals. {We remark that the weaker assumptions on $g$ imply the use of different techniques and strategy when controlling the term $(III)$ in \eqref{displayIII}. In particular, we need a higher integrability condition, which is reflected by suitable hypotheses on exponents contained at the beginning of Section \ref{sechi}}.\\
Next, taking advantage of the higher integrability of solutions to \eqref{minprob}, we will be able to prove the higher differentiability result. The proof of Theorem \ref{thmhd} is achieved by means of difference
quotient method, that is quite natural when trying to establish higher differentiabilty
results and local gradient estimates (see for instance \cite{marcellini1989}).
{
\begin{rmk}\label{ThmII}
We remark that assuming  a stronger Besov regularity on the coefficients of $\partial_zF(x,z)$, namely hypothesis (F4') below, we are able to prove that solutions to \eqref{minprob} inherit a Besov regularity of the type $B^\frac{\alpha}{2}_{2,t}$, with $t \ge 1 $. \\
More precisely, let $\alpha,\gamma \in (0,1)$ and let $\sigma$ be the exponent defined in \eqref{sigma}. Let $F(x,z)$ satisfy (F1) and (F2) for exponents $1<p<q$ verifying \eqref{gapqpx}. Moreover, assume that 
there exists a sequence of measurable non-negative functions $g_k \in L^{\sigma}_{\text{loc}}(\Omega)$ such that
$$\displaystyle\sum_{k=1}^{\infty} \Vert g_k \Vert^{t}_{L^{\sigma}(\Omega)} < \infty,$$
for some $t \ge 1$,
and at the same time
\begin{equation*}
    |\partial_{z}F(x,z)-\partial_{z}F(y, z)| \leq |x-y|^{\alpha} (g_k(x)+g_k(y)) (1+|z|^{q-1}) \eqno{\rm{{ (F4')}}}
\end{equation*}
for a.e. $x,y \in \Omega$ such that $2^{-k} \text{diam}(\Omega) \leq |x-y| < 2^{-k+1}\text{diam}(\Omega)$ and for every $z \in \mathbb{R}^{N \times n}$.
Let $u\in W^{1,p}(\Omega, \mathbb{R}^N)$ be a solution to the minimum problem \eqref{minprob}.
Assuming that
\eqref{lav} holds for a ball $B_r \Subset \Omega$ with $r \le 1$ and for some $H >0$,
then
\item[(i)] $Du\in L^{\tilde{q}}_{\textrm{loc}}(\Omega, \R^{N \times n}),$  for all $q \le \tilde{q}<p+\dfrac{\min \lbrace \alpha,2\gamma\rbrace}{\Theta(1-\gamma)}$;
\\ \item[(ii)] $V_p(Du)\in B^\frac{\alpha}{2}_{2,t,\textrm{loc}}(\Omega, \R^{N \times n})$, for all $t \ge 1$.\\
\noindent The higher integrability result goes exactly in the same way of that of Theorem \ref{thmhi}. On the other hand, for the proof of the last part we refer to \cite{eleuteri.passarelli,Gentile,grimaldi0,grimaldi.ipocoana}.
\end{rmk}}

The paper is organized as follows. After recalling some notation and preliminary results in Section \ref{secnot}, we prove our main results, Theorems \ref{thmhi} and \ref{thmhd} in Section \ref{profthm}.

\section{Notation and preliminary results}
\label{secnot}
For the rest of the paper, we denote by $C$, $c$ general positive constants. Different occurrences from line to line will be still denoted using the same letters. Relevant dependencies on parameters will be emphasized using parentheses or subscripts. 
We denote by $B(x,r)=B_{r}(x)= \{ y \in \mathbb{R}^{n} : |y-x | < r  \}$ the ball centered at $x$ of radius $r$. We shall omit the dependence on the center and on the radius when no confusion arises.

For the auxiliary function $V_{p}$, we recall the following estimate (see the proof of \cite[Lemma 8.3]{giusti}): 
\begin{lem}\label{D1}
Let $1<p<+\infty$. There exists a constant $c=c(n,p)>0$ such that
\begin{center}
$c^{-1}(\mu^{2}+|\xi|^{2}+|\eta|^{2})^{\frac{p-2}{2}} \leq \dfrac{|V_{p}(\xi)-V_{p}(\eta)|^{2}}{|\xi-\eta|^{2}} \leq c(\mu^{2}+|\xi|^{2}+|\eta|^{2})^{\frac{p-2}{2}} $
\end{center}
for any $\xi, \eta \in \mathbb{R}^{k}$, $\xi \neq \eta$, $k \in \N$.
\end{lem}

{Now we state a well-known iteration lemma (see \cite[Lemma 6.1]{giusti} for the proof).
\begin{lem}\label{lemreab}
Let $\Phi  :  [\frac{R}{2},R] \rightarrow \mathbb{R}$ be a bounded nonnegative function, where $R>0$. Assume that for all $\frac{R}{2} \leq r < s \leq R$ it holds
$$\Phi (r) \leq \theta \Phi(s) +A + \dfrac{B}{(s-r)^2}+ \dfrac{C}{(s-r)^{\gamma}}$$
where $\theta \in (0,1)$, $A$, $B$, $C \geq 0$ and $\gamma >0$ are constants. Then there exists a constant $c=c(\theta, \gamma)$ such that
$$\Phi \biggl(\dfrac{R}{2} \biggr) \leq c \biggl( A+ \dfrac{B}{R^2}+ \dfrac{C}{R^{\gamma}}  \biggr).$$
\end{lem}}

\subsection{Difference quotient}
\label{secquo}
We recall some properties of the finite difference quotient operator that will be needed in the sequel. For every function $F:\mathbb{R}^{n}\rightarrow \mathbb{R}^k$, $k \in \N$, let  
\begin{align*}
\tau_{h}F(x) :=F(x+h)-F(x),
\end{align*}
where $h \in \mathbb{R}^{n}$.
\\We start with the description of some elementary properties that can be found, for example, in \cite{giusti}.
\begin{prop}\label{rapportoincrementale}
Let $F \in W^{1,p}(\Omega, \R^k)$, with $p \geq1$, and let us consider the set
\begin{center}
$\Omega_{|h|} = \{ x \in \Omega : \mathrm{dist}(x,\partial \Omega)> |h|  \}$.
\end{center}
Then
\\(i) $\tau_{h}F \in W^{1,p}(\Omega_{|h|}, \R^k)$ and 
\begin{center}
$D_{i}(\tau_{h}F)=\tau_{h}(D_{i}F)$.
\end{center}
(ii) If at least one of the functions $F$ or $G$ has support contained in $\Omega_{|h|}$, then
\begin{center}
$$\displaystyle\int_{\Omega}F \tau_h G \ \mathrm {d}x = -\displaystyle\int_{\Omega} G \tau_{-h}F \ \mathrm{d}x.$$
\end{center}
(iii) We have $$\tau_h (FG)(x)= F(x+h)\tau_h G(x)+G(x) \tau_h F(x).$$
\end{prop}
The next result about finite difference operator is a kind of integral version of Lagrange Theorem.
\begin{lem}\label{ldiff}
If $0<\rho<R,$ $|h|<\frac{R-\rho}{2},$ $1<p<+\infty$ and $F \in W^{1,p}(B_{R}, \R^{k })$, then
\begin{center}
$\displaystyle\int_{B_{\rho}} |\tau_{h}F(x)|^{p} \ \mathrm{d} x \leq c(n,p)|h|^{p} \displaystyle\int_{B_{R}} |DF(x)|^{p} \ \mathrm{d}x$.
\end{center}
Moreover,
\begin{center}
$\displaystyle\int_{B_{\rho}} |F(x+h)|^{p} \ \mathrm{d} x \leq  \displaystyle\int_{B_{R}} |F(x)|^{p} \ \mathrm{d}x$.
\end{center}
\end{lem}

\subsection{Function spaces}
\label{secbesov}

Let $v:\mathbb{R}^{n} \rightarrow \mathbb{R}^k$ be a function, $k \in \N$. As in \cite[Section 2.2.2]{triebel83}, given $0<\alpha <1 $ and $1 \leq p,s< \infty$, we say that $v$ belongs to the Besov space $B^{\alpha}_{p,s}(\mathbb{R}^{n},\R^k)$ if $v \in L^{p}(\mathbb{R}^{n}, \R^k)$ and
\begin{center}
$[v]_{B^{\alpha}_{p,s}(\mathbb{R}^{n})} :=  \biggl( \displaystyle\int_{\mathbb{R}^{n}} \biggl( \displaystyle\int_{\mathbb{R}^{n}} \dfrac{|\tau_hv(x)|^{p}}{|h|^{\alpha p}} \ \mathrm{d} x \biggr)^{\frac{s}{p}}  \dfrac{  \mathrm{d} h}{|h|^{n}} \biggr)^{\frac{1}{s}}  < \infty$.
\end{center}
Equivalently, we could simply say that $v \in L^{p}(\mathbb{R}^{n},\R^k)$ and $\frac{\tau_{h}{v}}{|h|^{\alpha}} \in L^{s}\bigl( \frac{dh}{|h|^{n}}; L^{p}(\mathbb{R}^{n}) \bigr)$. 
Similarly, we say that $v \in B^{\alpha}_{p,\infty}(\mathbb{R}^{n},\R^k)$ if $v \in L^{p}(\mathbb{R}^{n},\R^k)$ and
\begin{center} 
$[v]_{B^{\alpha}_{p, \infty}(\mathbb{R}^{n})} : =  \displaystyle\sup_{h \in \mathbb{R}^{n}} \biggl( \displaystyle\int_{\mathbb{R}^{n}} \dfrac{|\tau_hv(x)|^{p}}{|h|^{\alpha p}} \ \mathrm{d} x \biggr)^{\frac{1}{p}} < \infty $.
\end{center}
We also define the following norm for $B^\alpha_{p,s}(\R^n,\R^k)$

\begin{equation}\label{besovnorm}
\Vert v \Vert_{B^{\alpha}_{p,s}(\mathbb{R}^{n})} : = \Vert v \Vert_{L^{p}(\mathbb{R}^{n})} + [v]_{B^{\alpha}_{p,s}(\mathbb{R}^{n})} .
\end{equation}
Moreover,
if one integrates for $h \in B(0, \delta)$ for a fixed $\delta >0$ then an equivalent norm is obtained, because
\begin{align}\label{equivnorm}
\biggl( \displaystyle\int_{\{|h| \geq \delta\}} \biggl( \displaystyle\int_{\mathbb{R}^{n}} \dfrac{|\tau_h v(x)|^{p}}{|h|^{\alpha p}} \ \mathrm{d}x \biggr)^{\frac{s}{p}}  \dfrac{ \mathrm{d}h}{|h|^{n}} \biggr)^{\frac{1}{s}} \leq c(n, \alpha,p,s, \delta) \Vert v \Vert_{L^{p}(\mathbb{R}^{n})} .
\end{align}
When $s= \infty$, one can simply take supremum over $|h| \leq \delta$ and obtain an equivalent norm.

 By construction, one has $B^{\alpha}_{p, s}(\mathbb{R}^{n},\R^k) \subset L^{p}(\mathbb{R}^{n},\R^k)$. One also has the following version of Sobolev embeddings (a proof can be found at \cite[Proposition 7.12]{haroske}).
\begin{lem}\label{3.1}
Suppose that $0 < \alpha <1$.
\\ (a) If $1 < p < \frac{n}{\alpha}$ and $1 \leq s \leq p^{*}_{\alpha} = \frac{np}{n- \alpha p}$, then there is a continuous embedding $B^{\alpha}_{p, s}(\mathbb{R}^{n}) \subset L^{p^{*}_{\alpha}}(\mathbb{R}^{n})$.
\\ (b) If $p = \frac{n}{\alpha}$ and $1 \leq s \leq \infty$, then there is a continuous embedding $B^{\alpha}_{p, s}(\mathbb{R}^{n}) \subset BMO(\mathbb{R}^{n})$,
\\ where $BMO$ denotes the space of functions with bounded mean oscillations \emph{\cite[Chapter 2]{giusti}}.
\end{lem}
We recall the following inclusions (\cite[Proposition 7.10 and Formula (7.35)]{haroske}).
\begin{lem}\label{3.2}
Suppose that $0 < \beta < \alpha <1$.
\\ (a) If $1 < p < \infty$ and $1 \leq s \leq t \leq \infty$, then $B^{\alpha}_{p, s}(\mathbb{R}^{n}) \subset B^{\alpha}_{p, t}(\mathbb{R}^{n})$.
\\ (b) If $1 < p < \infty$ and $1 \leq s , t \leq \infty$, then $B^{\alpha}_{p, s}(\mathbb{R}^{n}) \subset B^{\beta}_{p, t}(\mathbb{R}^{n})$.
\end{lem}

Given a domain $\Omega \subset \mathbb{R}^{n}$, we say that $v$ belongs to the local Besov space $ B^{\alpha}_{p, s,loc}$ if $\varphi \ v \in B^{\alpha}_{p, s}(\mathbb{R}^{n})$ whenever $\varphi \in \mathcal{C}^{\infty}_{c}(\Omega)$. It is worth noticing that one can prove suitable version of Lemmas \ref{3.1} and \ref{3.2}, by using local Besov spaces.

We have the following characterization for local Besov spaces (see \cite{baison.clop2017} for the proof).
\begin{lem}
A function $v \in L^{p}_{loc}(\Omega, \mathbb{R}^k)$ belongs to the local Besov space $B^{\alpha}_{p,s,loc}$, with $0< \alpha < 1$, if, and only if,
\begin{center}
$\biggl\Vert \dfrac{\tau_{h}v}{|h|^{\alpha}} \biggr\Vert_{L^{s}\bigl(\frac{dh}{|h|^{n}};L^{p}(B)\bigr)}<  \infty,$
\end{center}
for any ball $B\subset2B\subset\Omega$ with radius $r_{B}$, {where
\begin{align*}
\biggl\Vert \dfrac{\tau_{h}v}{|h|^{\alpha}} \biggr\Vert_{L^{s}\bigl(\frac{dh}{|h|^{n}};L^{p}(B)\bigr)}:=\biggl( \displaystyle\int_{\{|h| \geq \delta\}} \biggl( \displaystyle\int_{B} \dfrac{|\tau_h v(x)|^{p}}{|h|^{\alpha p}} \ \mathrm{d}x \biggr)^{\frac{s}{p}}  \dfrac{ \mathrm{d}h}{|h|^{n}} \biggr)^{\frac{1}{s}}. 
\end{align*} }
Here the measure $\frac{dh}{|h|^n}$ is restricted to the ball $B(0,r_B)$ on the h-space.
\end{lem}

It is known that Besov-Lipschitz spaces of fractional order $\alpha \in (0,1)$ can be characterized in pointwise terms. Adopting the terminology of \cite{koskela}, given a measurable function $v:\mathbb{R}^{n} \rightarrow \mathbb{R}$, a \textit{fractional $\alpha$-Hajlasz gradient for $v$} is a sequence $\{g_{k}\}_{k}$ of measurable, non-negative functions $g_{k}:\mathbb{R}^{n} \rightarrow \mathbb{R}$, together with a null set $N\subset\mathbb{R}^{n}$, such that the inequality 
\begin{center}
$|v(x)-v(y)|\leq (g_{k}(x)+g_{k}(y))|x-y|^{\alpha}$
\end{center} 
holds whenever $k \in \mathbb{Z}$ and $x,y \in \mathbb{R}^{n}\setminus N$ are such that $2^{-k} \leq|x-y|<2^{-k+1}$. We say that $\{g_{k}\}_{k} \in l^{s}(\mathbb{Z};L^{p}(\mathbb{R}^{n}))$ if
\begin{center}
$\Vert \{g_{k}\}_{k} \Vert_{l^{s}(L^{p})}=\biggl(\displaystyle\sum_{k \in \mathbb{Z}}\Vert g_{k} \Vert^{s}_{L^{p}(\mathbb{R}^{n})} \biggr)^{\frac{1}{s}}<  \infty.$
\end{center} 
In fact, $\{g_{k}\}_{k}$ above is not really a gradient. One should view it as a maximal function of the usual gradient.

The following result was proved in \cite{koskela}.
\begin{thm}\label{thm2.7}
Let $0< \alpha <1,$ $1 \leq p < \infty$ and $1\leq s \leq \infty $. Let $v \in L^{p}(\mathbb{R}^{n})$. One has $v \in B^{\alpha}_{p,s}(\mathbb{R}^{n})$ if, and only if, there exists a fractional $\alpha$-Hajlasz gradient $\{g_{k}\}_{k} \in l^{s}(\mathbb{Z};L^{p}(\mathbb{R}^{n}))$ for $v$. Moreover,
\begin{center}
$\Vert v \Vert_{B^{\alpha}_{p,s}(\mathbb{R}^{n})}\simeq \inf \Vert \{g_{k}\}_{k} \Vert_{l^{s}(L^{p})},$
\end{center}
where the infimum runs over all possible fractional $\alpha$-Hajlasz gradients for $v$.
\end{thm}

Let $\sigma \in (0, +\infty) \setminus \mathbb{N}$, $p \in [1, + \infty)$, $k \in \mathbb{N}$, $n \ge 2$. If $\sigma \in (0,1)$, we say that a function $v : \Omega \rightarrow \mathbb{R}^k$ belongs to the fractional Sobolev space $W^{\sigma, p}(\Omega, \mathbb{R}^k)$ if, and only if, the following Gagliardo type norm is finite:
\begin{align*}
\Vert v \Vert_{W^{\sigma, p}(\Omega)}:= &\Vert v \Vert_{L^p(\Omega)} + \biggl( \int_\Omega \int_\Omega \dfrac{|v(x)-v(y)|^p}{|x-y|^{n+\sigma  p}} \ \mathrm{d}x \mathrm{d}y  \biggr)^\frac{1}{p} \\
=: & \Vert v \Vert_{L^p(\Omega)} + [v]_{\sigma,p;\Omega}.
\end{align*}
On the other hand, if $\sigma = [\sigma]+ \{\sigma \} \in \mathbb{N}+ (0,1)>1$, we have $v \in W^{\sigma, p}(\Omega, \mathbb{R}^k)$ if, and only if, 
$$  \Vert v \Vert_{W^{\sigma, p}(\Omega)}:= \Vert v \Vert_{W^{[\sigma], p}(\Omega)}+ [D^{[\sigma]}v]_{\{\sigma \},p;\Omega}$$
is finite.

From Lemma \ref{3.2}, we deduce the following fundamental embedding property between Besov and fractional Sobolev spaces.
\begin{lem}\label{lem1DFM}
We consider radii $0<\rho<R$ and a function
  $v\in L^{p}(B_{R},\R^k)$.
  Suppose that there exist $\sigma\in (0,1)$
  and $M>0$ such that
  \begin{align*}
    \Vert \tau_h v \Vert_{L^p(B_\rho)} \le M |h|^{\sigma}
\end{align*}
for every $h\in\R^n$ with $|h|<\frac{R-\rho}{2}$. Then $v\in
W^{\beta,q}(B_{\rho}, \mathbb{R}^{k})$ for every $\beta\in (0,\sigma)$. Moreover,
\begin{equation*}
\Vert v\Vert_{W^{\beta,p}(B_\rho,\R^k)}
\leq \dfrac{c}{(R-\rho)^\delta} \bigl(
M+\Vert v \Vert_{L^p(B_R,\R^k)} \bigr),
\end{equation*}
holds true for constants $c,\delta$ depending on $n,p,\beta$ and $\sigma$.
\end{lem}

We conclude this section with a Gagliardo-Nirenberg type inequality in fractional Sobolev spaces 
(see \cite{DFM} for the proof).


\begin{lem}\label{lem3DFM} 
Let $B_\rho\Subset B_r\Subset \R^n$ be concentric balls with $r\leq 1$. Let $0\leq s_1<1<s_2<2, 1<a, t<\infty, \tilde{p}>1$ and $\theta\in(0,1)$ be such that
\begin{equation*}
\theta s_1+(1-\theta)s_2=1, \qquad \dfrac{1}{\tilde{p}}=\dfrac{\theta}{a}+\dfrac{1-\theta}{t}.
\end{equation*}
Then every function $v\in W^{s_1,a}(B_r,\R^k)\cap W^{s_2,t}(B_r,\R^k)$ belongs to $W^{1,\tilde{p}}(B_\rho,\R^N)$ and the inequality
\begin{equation*}
\Vert Dv\Vert_{L^{\tilde{p}}(B_\rho)}
\leq \dfrac{c}{(r-\rho)^\kappa}[v]^\theta_{s_1,a;B_r}
\Vert Dv\Vert_{W^{s_2-1,t}(B_r)}^{1-\theta}
\end{equation*}
holds for constants $c,\kappa=c,\kappa(n,s_1,s_2,a,t)$.
\end{lem}

\section{Proof of Theorems \ref{thmhi} and \ref{thmhd}}\label{profthm}
In the following we prove the main results of this work, namely the higher integrability result, Theorem \ref{thmhi}, and the higher differentiability one, Theorem \ref{thmhd}. The strategy proposed to prove Theorem \ref{thmhi} follows the one introduced in \cite{DFM}. However, we remark that we are here dealing with a weaker assumption on the summability of the map $x \mapsto\partial_z F(x,z) $ (see \eqref{sigma}). We denote by $u\in W^{1,p}(\Omega, \R^N) $ a solution to \eqref{minprob}. Throughout this section, we will assume that \eqref{lav} holds for a ball $B_r \Subset \Omega$ with $r \le 1$ and for some $H >0$.

\subsection{Higher integrability}\label{sechi}

\noindent We start by considering parameters $s,d$ and $\beta$ such that
\begin{equation}\label{param}
0\leq s<\gamma, \qquad 2d>\max\lbrace q,n\rbrace, \qquad 0<\beta<\min\lbrace \alpha,2\gamma\rbrace.
\end{equation}
We define the function $\tilde{p}=\tilde{p}(s,d,\beta)$ as
\begin{eqnarray}\label{ptil}
 \tilde{p}:=
  \begin{cases}
  \dfrac{2d(p(1-s)+\beta)}{\beta+2d(1-s)} & \textrm{if}\quad  p\geq 2  \bigskip\\
  \dfrac{2dp[2(1-s)+\beta]}{p\beta+4d(1-s)} & \textrm{if}\quad  1<p<2
  \end{cases} 
\end{eqnarray}

\noindent By definition, we have
\begin{equation}\label{dp}
d>\dfrac{p}{2},
\end{equation}
which yields
\begin{equation}\label{gapptild}
\tilde{p}<2d.
\end{equation}
The function $\tilde{p}$ is increasing in all its variables and
\begin{equation}
\lim_{s\to \gamma, d\to\infty, \beta\to \min\lbrace\alpha,2\gamma\rbrace}\tilde{p}(s,d,\beta)=p+\dfrac{\min\lbrace \alpha,2\gamma\rbrace}{\Theta(1-\gamma)}.
\end{equation}
We take $d$ large such that
\begin{equation}\label{beta0}
d>\dfrac{\max\lbrace q,n\rbrace}{2}\Rightarrow \beta_0:=s-\dfrac{n}{2d}>0.
\end{equation}
Now we choose $s,d$, satisfying \eqref{beta0}, and $\beta$ such that
\begin{equation}\label{beta}
\beta<\alpha_0:=\min\lbrace\alpha,2\beta_0\rbrace<\min\lbrace \alpha,2\gamma\rbrace
\end{equation}
and such that
\begin{equation}\label{ptilgap}
p+\dfrac{\xi}{1-\gamma}<\tilde{p}(s,d,\beta)<p+\dfrac{\min\lbrace \alpha,2\gamma\rbrace}{1-\gamma}.
\end{equation}

Now, we want to approximate, on a fixed ball $B_r \Subset \Omega$, the original minimizer $u$ of $\mathcal{F}$ with a sequence $\{ u_j \}$ of $W^{1,q}$-regular solutions to approximating variational problems. To find such $\{ u_j \}$, we define a sequence of suitable functionals $\mathcal{F}_j$ of mixed local/non-local type as follows. Consider the parameters $s$, $\gamma$ and $d$ defined above.  For every $j \in \mathbb{N}$, let $ \bar{v}_j \in W^{1,2d}(B_r,\R^N)\cap W^{s,2d}(\R^n,\R^N)$ and let introduce the non-local Dirichlet class
\begin{equation*}
\mathbb{X}(\bar{v}_j,B_r):=
\left\lbrace v\in \left( \bar{v}_j+W_0^{1,2d}(B_r,\mathbb{R}^N)\right) \cap W^{s,2d}(\R^n, \R^N): v\equiv \bar{v}_j \textrm{ on } \R^n\setminus B_r\right\rbrace.
\end{equation*}
This subset of $W^{1,2d}(B_r,\R^N)\cap W^{s,2d}(\R^n,\R^N)$ is convex and closed. Besides, it is non-empty, since
$\bar{v}_j \in \mathbb{X}(\bar{v}_j, B_r)$ by definition.
Then, we define 
\begin{align}\nonumber
\mathcal{F}_j(w, B_r):=&
\mathcal{F}(w, B_r)+\varepsilon_j \int_{B_r}(\mu^2+|Dw|^2)^d\ \mathrm{d}x\\
 &+
\int_{B_r}(|w|^2-M_0^2)^d_+ \ \mathrm{d}x
+\int_{\R^n}\int_{\R^n}
\dfrac{(|w(x)-w(y)|^2-M^2|x-y|^{2\gamma})^d_+}{|x-y|^{n+2sd}} \ \mathrm{d}x\mathrm{d}y,
\label{appfunctional}
\end{align}
for every $w \in \mathbb{X}(\bar{v}_j,B_r)$ and 
for some positive parameters $\varepsilon_j$, $M_0$
and $M$. 
 Thanks to assumption \eqref{lav}, that ensures the absence of the Lavrentiev phenomenon, we are able to prove the existence of good boundary values for which 
 the minimizers $ u_j's$ of \eqref{appfunctional} converge to $u$ and are bounded in the fractional Sobolev space $W^{s,2d}$, exactly as $u$. Indeed, the following result holds (for the proof we refer to  \cite[Sections 3.2 and 3.4]{DFM}).

 \begin{prop}\label{boundaryvalue}
 There exists a sequence $\{ \bar{v}_j \} \in W^{1,2d}(B_r,\R^N)\cap W^{s,2d}(\R^n,\R^N)$ such that  $[\bar{v}_j]_{0,\gamma; \mathbb{R}^n }\le c_*H$, with a constant $c_*:=c_*(n,d,\gamma,s)$.
 Moreover, choosing in \eqref{appfunctional} parameters
\begin{equation}
\begin{cases}
 \varepsilon_j:=
 \dfrac{1}{\left(\Vert D\bar{v}_j\Vert^{4d}_{L^{2d}(B_r)}+j+1\right)} \smallskip  \\
   M_0:= 16c_*r^\gamma H \smallskip\\
    M:= 16c_* H
\end{cases}
\end{equation}
we have that, for every $j \in \N$, there exists $u_j \in \mathbb{X}(\bar{v}_j, B_r)$  solution
to 
\begin{equation*}\label{solution}
{\mathcal{F}_j(u_j, B_r)\leq \mathcal{F}_j(\omega, B_r) \qquad \textrm{for all } \omega \in \mathbb{X}(\bar{v}_j, B_r)}
\end{equation*} 
such that
\begin{equation}\label{convtou}
    u_j \rightharpoonup u \quad \text{weakly in} \quad W^{1,p} \bigl( B_r, \mathbb{R}^N \bigr)
\end{equation}
and for every $j \in \N$ it holds
\begin{equation}\label{holdernorm}
     [u_j]^{2d}_{s,2d;\R^n} \le c r^{n+2d (\gamma-s)} H^{2d} ,
\end{equation}
for a constant $c:=c(n,p,q,\alpha,\gamma, \sigma)$.\\
Moreover,
\begin{equation}\label{estFj}
\varepsilon_j\int_{B_r}\left(\mu^2+|Du_j|^2\right)^d \mathrm{d} x
\leq \mathcal{F}(u, B_r)+1
\end{equation}
and
\begin{equation}\label{estujFmat}
\Vert u_j\Vert_{W^{1,p}(B_r)}\leq c[\mathbb{F}(u,B_r)]^\frac{1}{p},
\end{equation}
where 
\begin{equation}\label{defFmat}
\mathbb{F}(u,B_r): =
\mathcal{F}(u,B_r)+r^nH^{2d}+1,
\end{equation}
for $c:=c(n,N,p,q,\nu,L,\alpha,\gamma, \sigma)$.
\end{prop}

We are now in position to give the proof of the higher integrability result.
\\

\noindent \textit{Proof of Theorem \ref{thmhi}.} For every $j \in \N$, let $u_j \in \mathbb{X}(\bar{v}_j, B_r)$ be the minimizer of \eqref{appfunctional}, where $\{ \bar{v}_j \} \in W^{1,2d}(B_r,\R^N)\cap W^{s,2d}(\R^n,\R^N)$ and parameters $\varepsilon_j$, $M_0$ and $M$ are as in Proposition \ref{boundaryvalue}.
We derive the Euler-Lagrange system of the functional $\mathcal{F}_j$ defined in \eqref{appfunctional}.  For every $j \in \mathbb{N}$ and $t \in (-1,1)$, we consider the variation $u_j+t \varphi$, where $\varphi \in W^{1,2d}_0(B_r,\R^N)$ has compact support in $B_r$. Note that, from the embedding $W^{1,2d}(B_r,\R^N) \subset W^{s,2d}(B_r,\R^N)$, we have $\varphi \in  W^{s,2d}(B_r,\R^N)$. Moreover, since $\varphi$ has compact support in $B_r$, setting $\varphi=0$ outside $B_r$ yields that $\varphi \in  W^{s,2d}(\R^n,\R^N)$. Therefore, we have that
$\varphi \in W^{1,2d}_0(B_r,\R^N)\cap W^{s,2d}_0(\R^n,\R^N)$ and $u_j+t \varphi \in \mathbb{X}(\bar{v}_j,B_r)$. By minimality of $u_j$, we get 
$$\dfrac{\mathrm{d} \mathcal{F}_j(u_j+t \varphi, B_r)}{\mathrm{d}t} \bigg|_{t=0}=0,$$
that is
\begin{align}\nonumber
0 =&
\int_{B_r}{\partial_z}F_j(x, Du_j) \cdot D \varphi \ \mathrm{d}x+2d \int_{B_r}(|u_j|^2-M_0^2)_+^{d-1}\ u_j \cdot \varphi \ \mathrm{d}x\\
&+\int_{\R^n}\int_{\R^n}
\dfrac{(|u_j(x)-u_j(y)|^2-M^2|x-y|^{2\gamma})^{d-1}_+}{|x-y|^{n+2sd}}
(u_j(x)-u_j(y))\cdot (\varphi(x)-\varphi(y))
\ \mathrm{d}x\mathrm{d}y, \label{eulerlagrange}
\end{align}
where we used the notation
\begin{equation*}
F_j(x,z):=F(x,z)+\varepsilon_j(\mu^2+z^2)^d. \label{defFj}
\end{equation*}
Now, we fix radii
\begin{equation*}
    0< \varrho \le \tau_1 < \tau_2 \le r
\end{equation*}
and set 
\begin{equation*}
    \varphi:= \tau_{-h}(\eta^2\tau_hu_j) 
\end{equation*}
with $\eta \in \mathcal{C}^1_c(B_{(3\tau_2+\tau_1)/4}) $ verifying $1_{B_{(\tau_2+\tau_1)/2}} \le \eta \le 1_{B_{(3\tau_2+\tau_1)/4}}$ and $|D\eta| \le \frac{c}{\tau_2-\tau_1}$, and $h \in \mathbb{R}^n \setminus \{0\}$ is such that 
\begin{equation*}
    0 < |h| < \dfrac{\tau_2-\tau_1}{2^{10}}
\end{equation*}
Testing \eqref{eulerlagrange} with $\varphi$ and using integration by parts for finite difference operators,
we obtain
\begin{small}
\begin{align}\nonumber
\int_{\R^n} & \int_{\R^n}
\dfrac{\tau_h[(|u_j(x)-u_j(y)|^2-M^2|x-y|^{2\gamma})^{d-1}_+ (u_j(x)-u_j(y))]}{|x-y|^{n+2sd}}
\cdot (\eta^2(x)\tau_hu_j(x)-\eta^2(y)\tau_hu_j(y))
\mathrm{d}x\mathrm{d}y \\ \nonumber
+& 2d \int_{B_r}\eta^2 \tau_h \bigl((|u_j|^2-M_0^2)_+^{d-1}\ u_j \bigr) \cdot \tau_h u_j \mathrm{d}x \\ \label{sumI-II-III0}
+& \int_{B_r}\tau_h \partial_z F_j(x, Du_j) \cdot [\eta^2 \tau_h D u_j + 2 \eta \tau_h u_j \otimes D \eta] \mathrm{d}x := (I)+ (II)+(III)=0.
\end{align}
\end{small}

\noindent Concerning the term $(I)$, we take advantage of estimate in \cite[Section 3.6]{DFM}, in particular it holds

\begin{equation}\label{(I)}
d\int_{\mathbb{R}^n}\int_{\mathbb{R}^n}\int_0^1 \dfrac{\eta^2(x)B_1(x,y,\lambda)}{|x-y|^{n+2sd}}\dd \lambda \dd x \dd y\leq (I)+\dfrac{cr^{\frac{n}{d}}\mathbb{F}(u,B_r)}{(\tau_2-\tau_1)^2}|h|^{2\beta_0},
\end{equation}
where we defined
$$B_1(x,y,\lambda):= \biggl[\dfrac{d}{d \lambda} (|U_{\lambda,h}(x,y)|^2-M^2|x-y|^{2 \gamma})^{d-1}_+ U_{\lambda,h}(x,y) \biggr] \cdot (\tau_hu_j(x)-\tau_hu_j(y)) ,$$
with 
$$U_{\lambda,h}(x,y):= u_j(x)-u_j(y)+\lambda \tau_hu_j(x)-\lambda \tau_hu_j(y)$$
and where $\mathbb{F}$ was introduced in \eqref{defFmat} and $\beta_0= s-n/2d $.
\\Next, we estimate $(II)$ as follows
\begin{align}\nonumber
(II)=& \ 2d  \int_{B_r}\eta^2 \tau_h u_j\int_0^1 \dfrac{\mathrm{d}}{\mathrm{d}\lambda} \biggl( 
\bigl(  |u_j+\lambda \tau_h u_j|^2-M_0^2 \bigr)^{d-1}_+(u_j+\lambda \tau_hu_j)
\biggr)\ \mathrm{d}\lambda \mathrm{d}x
\\
=&\ 4d(d-1)\int_{B_r}\eta^2\int_0^1\left(|u_j+\lambda\tau_hu_j|^2-M^2\right)^{d-2}_+((u_j+\lambda\tau_h u_j)\cdot \tau_h u_j)^2\dd \lambda \dd x \notag\\ \label{II}
&+ 2d\int_{B_r}\eta^2|\tau_h u_j|^2\int_0^1\left(|u_j+\lambda\tau_hu_j|^2-M_0^2\right)^{d-1}_+\dd \lambda\dd x \geq 0.
\end{align}

Now we take care of the integral $(III)$. Recalling the definition \eqref{defFj}, we have 

\begin{align}\nonumber
(III) = &\int_{B_r}\tau_h \partial_z F(x,Du_j) \cdot [\eta^2\tau_h Du_j+2\eta\tau_h u_j\otimes D\eta] \ \mathrm{d}x\\ \label{displayIII}
&+ 2d \varepsilon_j \int_{B_r} ((\mu^2+|Du_j|^2)^{d-1}Du_j) \cdot [\eta^2\tau_h Du_j+2\eta\tau_h u_j\otimes D\eta] \ \mathrm{d}x =:(III)_1+(III)_2.
\end{align}
For $(III)_1$, we decompose
\begin{align}\nonumber
(III)_1=&\int_{B_r}[\partial_z F(x+h,Du_j(x+h))-\partial_z F(x+h,Du_j(x))] \cdot (\eta^2\tau_h Du_j+2\eta {\tau_h u_j\otimes D\eta}) \ \mathrm{d}x\\ \nonumber
&+\int_{B_r}[\partial_z F(x+h,Du_j(x))-\partial_z F(x,Du_j(x))] \cdot (2\eta {\tau_h u_j\otimes D\eta} +\eta^2\tau_h Du_j) \ \mathrm{d}x\\ \label{III1}
=&:(III)_{1,1}+(III)_{1,2}+(III)_{1,3}+(III)_{1,4}
\end{align}
By the ellipticity assumption (F2), we have
\begin{equation}\label{estIII11}
(III)_{1,1}\geq \dfrac{1}{c}\int_{B_r}\eta^2(\mu^2+|Du_j(x+h)|^2+|Du_j(x)|^2)^\frac{p-2}{2} |\tau_h Du_j|^2 \ \mathrm{d}x.
\end{equation}
Using growth assumption (F1) and then H\"{o}lder's inequality, we obtain
 \begin{align}\nonumber
 (III)_{1,2}+(III)_{1,3}&\leq c \int_{B_r} \eta (1+|Du_j(x+h)|^2+|Du_j(x)|^2)^\frac{q-1}{2}|D\eta||\tau_h u_j| \ \mathrm{d}x\\ \nonumber
 &\leq c \Vert D\eta \Vert_{L^\infty}\left(\int_{B_r}\eta(1+|Du_j(x+h)|^2+|Du_j(x)|^2)^\frac{q}{2} \ \mathrm{d}x\right)^\frac{q-1}{q}\\ \nonumber
 &\qquad \cdot\left(\ir \eta|\tau_h u_j|^q \ \mathrm{d} x\right)^\frac{1}{q}\\ \nonumber
  &\leq \dfrac{c}{\tau_2-\tau_1}|h|\it2(1+|Du_j|^2)^\frac{q}{2}\ \mathrm{d}x\\ \label{estIII123}
   &\leq \dfrac{c}{\tau_2-\tau_1}|h|\it2(1+|Du_j|)^{p+\frac{\xi}{1-\gamma}}\ \mathrm{d}x,
 \end{align}
where $\xi$ was {introduced} in \eqref{sigma} and we have also exploited the bounds { $q<p+\frac{\xi}{1-\gamma}$} and \eqref{gapqpx} and hypotheses on $\eta$ and $D\eta$.\\
According to (F3) we can estimate $(III)_{1,4}$ as
$$(III)_{1,4} \leq c |h|^\alpha \ir \eta^2(g(x+h)+g(x))(1+|Du_j|^2)^\frac{q}{2}\ \mathrm{d}x.$$
We now use H\"{o}lder's inequality with exponents $\sigma$ and $\dfrac{\sigma}{\sigma-1}$, where $\sigma$ was introduced in \eqref{sigma}. Then, by \eqref{gapqpx} and Young's inequality we infer
\begin{align}\label{estIII14}
(III)_{1,4}\leq |h|^\alpha \left(\ir g^\sigma \ \mathrm{d}x\right)^\frac{1}{\sigma}
\left(\it2 ( 1+|Du_j|^{p+\frac{\xi}{1-\gamma}}) \ \mathrm{d}x\right)^{\frac{q(1-\gamma)}{p(1-\gamma)+\xi}}.
\end{align}
{We notice that in the last in integral in the right hand side, we exploited the fact that the choices of the parameters
ensure $\frac{\sigma}{\sigma-1}q= p+\frac{\xi}{1-\gamma}$.
}
Putting \eqref{estIII11}, \eqref{estIII123} and \eqref{estIII14} in \eqref{III1}, we get
\begin{align}\nonumber
(III)_1\geq &\ \dfrac{1}{c}\int_{B_r}\eta^2(\mu^2+|Du_j(x+h)|^2+|Du_j(x)|^2)^\frac{p-2}{2} |\tau_h Du_j|^2 \ \mathrm{d}x\\ \label{estIII1}
 &- c\ir g^\sigma dx - \dfrac{c}{\tau_2-\tau_1}|h|^\alpha\it2(1+|Du_j|)^{p+\frac{\xi}{1-\gamma}} \ \mathrm{d}x.
\end{align}
On the other hand, similarly we might estimate
\begin{align}\nonumber
(III)_2\geq & \ \dfrac{\varepsilon_j}{c}\int_{B_r}\eta^2\left(\mu^2+|Du_j(x+h)|^2+|Du_j(x)|^2\right)^{d-1}|\tau_hDu_j|^2\dd x\\ \label{estIII2}
& \qquad -\dfrac{c\varepsilon_j|h|}{\tau_2-\tau_1}\int_{B_{\tau_2}}\left(1+|Du_j|^2\right)^d \dd x,
\end{align}
where we recall that $\varepsilon_j$ was defined in Proposition \ref{boundaryvalue}.
Using the estimates for $(I), (II)$ and $(III)$ in \eqref{sumI-II-III0} and exploiting \eqref{beta} and the hypotheses on $D\eta$ together with \eqref{estFj} and \eqref{defFmat}, we conclude
\begin{align}\nonumber
&\int_{B_\frac{\tau_2+\tau_1}{2}}\int_{B_\frac{\tau_2+\tau_1}{2}}\int_0^1 \dfrac{B_1(x,y,\lambda)}{|x-y|^{n+2sd}} \ \mathrm{d}\lambda  \mathrm{d}x \mathrm{d}y\\ \nonumber
&\qquad +\int_{B_\frac{\tau_2+\tau_1}{2}}(\mu^2+|Du_j(x+h)|^2+|Du_j(x)|^2)^\frac{p-2}{2} |\tau_h Du_j|^2 \ \mathrm{d}x\\ \nonumber
& \leq \dfrac{c}{(\tau_2-\tau_1)^2}|h|^{\alpha_0}\left( \mathbb{F}(u,B_r)+\Vert g\Vert^\sigma_{L^\sigma(B_r)}+\Vert Du_j\Vert^{p+\frac{\xi}{1-\gamma}}_{L^{p+\frac{\xi}{1-\gamma}}(B_{\tau_2})}+1\right),\\ \label{MI}
\end{align}
for a constant $c:=c(n,N,\nu,L,p,q)$ {and where $\alpha_0$ is defined in \eqref{beta}}.
From the definition of $B_1$, we can easily derive that $B_1(x,y,\lambda) \ge 0$. 
Hence, we can discard the first integral on the left hand side of \eqref{MI} and obtain
\begin{align}
& \int_{B_\frac{\tau_2+\tau_1}{2}} (\mu^2+|Du_j(x+h)|^2+|Du_j(x)|^2)^\frac{p-2}{2} |\tau_h Du_j|^2 \ \mathrm{d}x \nonumber\\
& \leq  \dfrac{c}{(\tau_2-\tau_1)^2}|h|^{\alpha_0}\left( \mathbb{\tilde{F}}(u,B_r)+\Vert Du_j\Vert^{p+\frac{\xi}{1-\gamma}}_{L^{p+\frac{\xi}{1-\gamma}}(B_{\tau_2})}+1\right), \label{MI1} 
\end{align}
where we set $\mathbb{\tilde{F}}(u,B_r) :=\mathbb{F}(u,B_r)+\Vert g\Vert^\sigma_{L^\sigma(B_r)}$.\\

Now, we distinguish the case $p \ge 2$ and the case $1 < p <2$.\\
Let $p\geq 2$. Then, from \eqref{MI1} we obtain 
\begin{align}\label{regularity}
\int_{B_\frac{\tau_2+\tau_1}{2}} |\tau_h Du_j|^p\leq\dfrac{c}{(\tau_2-\tau_1)^2}|h|^{\alpha_0}\left( \mathbb{\tilde{F}}(u,B_r)+\Vert Du_j\Vert^{p+\frac{\xi}{1-\gamma}}_{L^{p+\frac{\xi}{1-\gamma}}(B_{\tau_2})}+1\right).
\end{align}
From Lemma \ref{lem1DFM}, since $0<|h|<\dfrac{\tau_2-\tau_1}{2^{10}}$, according to \eqref{regularity} we read the following regularity of $u_j$
\begin{equation*}
u_j\in W^{1+\frac{\beta}{p},p}(B_\frac{\tau_2+\tau_1}{2},\R^n)\cap W^{s,2d}(B_r,\R^n),
\end{equation*}
for $\beta \in (0, \alpha_0)$ satisfying \eqref{ptilgap}.
Hence, by Lemma \ref{lem1DFM} and inequality \eqref{estujFmat}, we infer

\begin{align}\nonumber
\Vert u_j\Vert_{W^{1+\frac{\beta}{p},p}\left(B_\frac{\tau_2+\tau_1}{2}\right)} 
&\leq \dfrac{c}{(\tau_2-\tau_1)^{\delta}}
{\Vert u_j\Vert}_{W^{1,p}(B_{\tau_2})}\\ \nonumber
&\qquad +\dfrac{c}{(\tau_2-\tau_1)^{\frac{2}{p}+\delta}}\left( [\mathbb{\tilde{F}}(u,B_r)]^\frac{1}{p}+\Vert Du_j\Vert^{{1+\frac{\xi}{p-p\gamma}}}_{L^{p+\frac{\xi}{1-\gamma}}(B_{\tau_2})}+1\right)\\ \label{pma2-1}
&\leq \dfrac{c}{(\tau_2-\tau_1)^{\frac{2}{p}+\delta}}\left( [\mathbb{\tilde{F}}(u,B_r)]^\frac{1}{p}+\Vert Du_j\Vert^{{1+\frac{\xi}{p-p\gamma}}}_{L^{p+\frac{\xi}{1-\gamma}}(B_{\tau_2})}+1\right),
\end{align}
with $c,\delta=c,\delta(n,p,\beta,\alpha_0)$.
Now, applying Lemma \ref{lem3DFM}, we get
\begin{equation}\label{interplem3DFM}
\Vert Du_j\Vert_{L^{\tilde{p}}(B_{\tau_1})}\leq 
\dfrac{c}{(\tau_2-\tau_1)^\kappa} [u_j]_{s,2d,B_{\tau_2}}^\theta 
\Vert Du_j\Vert^{1-\theta}_{W^{\frac{\beta}{p},p}\left(B_\frac{\tau_2+\tau_1}{2}\right)}, 
\end{equation}
for $\theta=\dfrac{\beta}{p(1-s)+\beta}\in (0,1)$, where $c$ and $\kappa$ are the constants from Lemma \ref{lem3DFM}.
Moreover, we observe that from \eqref{holdernorm} it follows
\begin{equation}\label{st3.40DFM}
[u_j]_{s,2d,\R^n}\leq c[\mathbb{F}(u,B_r)]^\frac{1}{2d}.
\end{equation}

\noindent Therefore, combining \eqref{pma2-1}, \eqref{interplem3DFM} and \eqref{st3.40DFM}, we obtain

\begin{align} \nonumber
\Vert Du_j\Vert_{L^{\tilde{p}}(B_{\tau_1})}
&\leq 
\dfrac{c}{(\tau_2-\tau_1)^{\frac{2}{p}+\delta+\kappa}} [\mathbb{\tilde{F}}(u,B_r)]^\frac{1-\theta}{p}[\mathbb{F}(u,B_r)]^\frac{\theta}{2d}\\ \nonumber
&\qquad+ \dfrac{c}{(\tau_2-\tau_1)^{\frac{2}{p}+\delta+\kappa}} 
[\mathbb{F}(u,B_r)]^\frac{\theta}{2d}
\Vert Du_j\Vert^{\left(p+\frac{\xi}{1-\gamma}\right)\frac{1-\theta}{p}}_{L^{p+\frac{\xi}{1-\gamma}}(B_{\tau_2})}\\ \label{st-4}
&\qquad{+ \dfrac{c}{(\tau_2-\tau_1)^{\frac{2}{p}+\delta+\kappa}} 
[\mathbb{F}(u,B_r)]^\frac{\theta}{2d}}.
\end{align}

\noindent We now recall that by assumption (see \eqref{ptil} and \eqref{ptilgap})
\begin{equation}\label{controlpptil}
p<p+\frac{\xi}{1-\gamma}<\tilde{p}.
\end{equation}
Therefore, we can use the following interpolation inequality
\begin{equation}\label{interp-5}
\Vert Du_j\Vert_{L^{p+\frac{\xi}{1-\gamma}}}
\leq \Vert Du_j\Vert^{\theta_1}_{L^{\tilde{p}}}
\Vert Du_j\Vert^{1-\theta_1}_{L^p},
\end{equation}
where $\theta_1\in (0,1)$ is such that
\begin{equation*}
\dfrac{1-\gamma}{p(1-\gamma)+\xi}=\dfrac{\theta_1}{\tilde{p}}+\dfrac{1-\theta_1}{p},
\end{equation*}
hence
\begin{equation}\label{teta1}
\theta_1=\dfrac{\xi\tilde{p}}{(\tilde{p}-p)(p(1-\gamma)+\xi)}.
\end{equation}
Inserting \eqref{interp-5} in \eqref{st-4}, we infer
\begin{align*} \nonumber
\Vert Du_j\Vert_{L^{\tilde{p}}(B_{\tau_1})}
&\leq 
c [\mathbb{\tilde{F}}(u,B_r)]^\frac{1-\theta}{p}[\mathbb{F}(u,B_r)]^\frac{\theta}{2d}\\ \label{st4-5}
&\qquad+c
[\mathbb{F}(u,B_r)]^\frac{\theta}{2d}
\Vert Du_j\Vert^{\theta_1\left(p+\frac{\xi}{1-\gamma}\right)\frac{1-\theta}{p}}_{L^{\tilde{p}}(B_{\tau_2})}
\Vert Du_j\Vert^{(1-\theta_1)\left(p+\frac{\xi}{1-\gamma}\right)\frac{1-\theta}{p}}_{L^{p}(B_{\tau_2})}\\
& \qquad {+ c [\mathbb{F}(u,B_r)]^\frac{\theta}{2d}}
\end{align*}
Recalling the definition of $\tilde{p}$  in \eqref{ptil}, 
we note that when $p\geq 2$, it holds
\begin{align*}\nonumber
&\lim_{s\to \gamma, d\to\infty, \beta\to \min\lbrace\alpha,2\gamma\rbrace}
\dfrac{\tilde{p}\xi}{\tilde{p}-p}\dfrac{1}{1-\gamma}\dfrac{1-s}{p(1-s)+\beta}\\ \nonumber
\stackrel{\eqref{ptil}}{=}&
\lim_{s\to \gamma, d\to\infty, \beta\to \min\lbrace\alpha,2\gamma\rbrace}
\dfrac{2d[p(1-s)+\beta]}{\beta+2d(1-s)} \dfrac{\beta+2d(1-s)}{2d[p(1-s)+\beta]-p\beta-2pd(1-s)}\\ \nonumber
&\qquad \cdot \dfrac{\xi(1-s)}{(1-\gamma)[p(1-s)+\beta]}\\
&\stackrel{}{=}\lim_{s\to \gamma, d\to\infty, \beta\to \min\lbrace\alpha,2\gamma\rbrace}\dfrac{2d\xi}{1-\gamma}\dfrac{1-s}{\beta(2d-p)}=\dfrac{\xi}{\min\lbrace \alpha,2\gamma\rbrace}<1,
\end{align*}
where we remark that $\frac{\xi}{\min\lbrace \alpha,2\gamma\rbrace}$ can be estimated by $1$ according to (F3).
Therefore, we increase $s,d,\beta$ in order to have
\begin{equation}\label{incrsdbet}
\dfrac{\tilde{p}\xi}{\tilde{p}-p}\dfrac{1}{1-\gamma}\dfrac{1-s}{p(1-s)+\beta}<1.
\end{equation}
Now, owing to \eqref{teta} and \eqref{teta1},
\begin{align*}
\theta_1\left(p+\frac{\xi}{1-\gamma}\right)\frac{1-\theta}{p} = 
\dfrac{\tilde{p}\xi}{\tilde{p}-p}\dfrac{1}{1-\gamma}\dfrac{1-s}{p(1-s)+\beta}<1
\end{align*}
where we exploited \eqref{incrsdbet}. Applying Young's inequality and reabsorbing the term $\Vert Du_j\Vert_{L^{\tilde{p}}}$ {according to Lemma \ref{lemreab}}, we get
\begin{equation*}
\Vert Du_j\Vert_{L^{\tilde{p}}(B_{\tau_1})}\leq \dfrac{c}{(\tau_2-\tau_1)^{\kappa_1}} [\tilde{\mathbb{F}}(u,B_r)]^{\kappa_2},
\end{equation*}
for some exponents $\kappa_1$ and $\kappa_2$.
From \eqref{convtou} and weak lower semicontinuity, letting  $j\to +\infty$ we get 
$$Du\in L^{\tilde{p}}(B_{\tau_1}),$$
hence
\begin{equation*}
Du \in L^{\tilde{q}}_{\textrm{loc}}(\Omega),
\end{equation*}
for every $q<\tilde{q}<p+\dfrac{\min\lbrace \alpha,2\gamma\rbrace}{1-\gamma}$. This proves the theorem in the case $p \ge 2$.
For the higher integrability in the case
$1<p<2$ we proceed analogously.
Being $\tilde{p}$  defined as in \eqref{ptil}, for 
$1<p<2$ it holds

\begin{align}\nonumber
&\lim_{s\to \gamma, d\to\infty, \beta\to \min\lbrace\alpha,2\gamma\rbrace}
\dfrac{\tilde{p}\xi}{\tilde{p}-p}\dfrac{1}{1-\gamma}\dfrac{2(1-s)}{p[2(1-s)+\beta]}\\ \nonumber
=&
\lim_{s\to \gamma, d\to\infty, \beta\to \min\lbrace\alpha,2\gamma\rbrace}
\dfrac{\xi}{1-\gamma}\dfrac{4d(1-s)}{p\beta(2d-p)}=\dfrac{2\xi}{\min\lbrace \alpha,2\gamma\rbrace}<1,
\end{align}
where we used the assumption $\xi<\dfrac{\min\lbrace \alpha,2\gamma\rbrace}{\Theta}$ and the definition of $\Theta$ in \eqref{teta} for $1<p<2$.\\
Hence, we can find $s,d,\beta$ sufficiently large, such that
\begin{align}\label{estlt1p12}
\dfrac{\tilde{p}\xi}{(\tilde{p}-p)(1-\gamma)}\dfrac{2(1-s)}{p[2(1-s)+\beta]}<1.
\end{align}
Now, we note that
\begin{align*}\nonumber
\int_{B_\frac{\tau_2+\tau_1}{2}} |\tau_h Du_j|^p \ \mathrm{d}x 
&= \int_{B_\frac{\tau_2+\tau_1}{2}}
|\tau_h Du_j|^p(\mu^2+|Du_j(x+h)|^2+|Du_j(x)|^2)^\frac{(p-2)p}{4} \  \\&\qquad\qquad\qquad\cdot(\mu^2+|Du_j(x+h)|^2+|Du_j(x)|^2)^\frac{(2-p)p}{4} \ \mathrm{d}x.
\end{align*}
Using H\"{o}lder's inequality we derive
\begin{align*}\nonumber
\int_{B_\frac{\tau_2+\tau_1}{2}} |\tau_h Du_j|^p \ \mathrm{d}x 
&\leq\left(\int_{B_\frac{\tau_2+\tau_1}{2}}(\mu^2+|Du_j(x+h)|^2+|Du_j(x)|^2)^\frac{p-2}{2}
|\tau_h Du_j|^2 \ \mathrm{d}x\right)^\frac{p}{2}\\
&\qquad \cdot\left(\int_{B_\frac{\tau_2+\tau_1}{2}}
(\mu^2+|Du_j(x+h)|^2+|Du_j(x)|^2)^\frac{p}{2} \ \mathrm{d}x \right)^{1-\frac{p}{2}}.
\end{align*}
Inserting estimates \eqref{estujFmat} and \eqref{MI1} above, we obtain
\begin{align*}\nonumber
\int_{B_\frac{\tau_2+\tau_1}{2}} |\tau_h Du_j|^p \ \mathrm{d}x
&\leq  \left(\int_{B_\frac{\tau_2+\tau_1}{2}}(\mu^2+|Du_j(x+h)|^2+|Du_j(x)|^2)^\frac{p}{2}
 \ \mathrm{d}x\right)^\frac{p}{2} [\mathbb{F}(u,B_r)]^{1-\frac{p}{2}} \\
 &\leq \dfrac{c}{(\tau_2-\tau_1)^p}|h|^\frac{\alpha_0 p}{2}\left( \mathbb{\tilde{F}}(u,B_r)+\Vert Du_j\Vert^{p+\frac{\xi}{1-\gamma}}_{L^{p+\frac{\xi}{1-\gamma}}(B_{\tau_2})}+1\right).
\end{align*}
Hence, from Lemma \ref{lem1DFM} and $0<|h|<\dfrac{\tau_2-\tau_1}{2^{10}}$, we infer
\begin{equation*}
u_j\in W^{1+\frac{\beta}{2},p}(B_\frac{\tau_2+\tau_1}{2},\R^n)\cap W^{s,2d}(B_r,\R^n),
\end{equation*}
for $\beta \in (0, \alpha_0)$ satisfying \eqref{ptilgap} and \eqref{estlt1p12}.
Therefore, from Lemma \ref{lem3DFM}, it follows

\begin{align*}\nonumber
\Vert u_j\Vert_{W^{1+\frac{\beta}{p},p}\left(B_\frac{\tau_2+\tau_1}{2}\right)} 
&\leq \dfrac{c}{(\tau_2-\tau_1)^{\delta}}
{\Vert u_j\Vert}_{W^{1,p}(B_{\tau_2})}\\ \nonumber
&\qquad +\dfrac{1}{(\tau_2-\tau_1)^{1+\delta}}\left( [\mathbb{\tilde{F}}(u,B_r)]^\frac{1}{p}+\Vert Du_j\Vert^{p+\frac{\xi}{1-\gamma}}_{L^{p+\frac{\xi}{1-\gamma}}(B_{\tau_2})}+1\right)\\ \label{pma2-1}
&\leq \dfrac{1}{(\tau_2-\tau_1)^{1+\delta}}\left( [\mathbb{\tilde{F}}(u,B_r)]^\frac{1}{p}+\Vert Du_j\Vert^{p+\frac{\xi}{1-\gamma}}_{L^{p+\frac{\xi}{1-\gamma}}(B_{\tau_2})}+1\right).
\end{align*}
Now, applying Lemma \ref{lem3DFM}, we get
\begin{equation*}\label{interplem3DFMp12}
\Vert Du_j\Vert_{L^{\tilde{p}}(B_{\tau_1})}\leq 
\dfrac{c}{(\tau_2-\tau_1)^\kappa} [u_j]_{s,2d,B_{\tau_2}}^\theta 
\Vert u_j\Vert^{1-\theta}_{W^{\frac{\beta}{p},p} \bigl(B_\frac{\tau_2+\tau_1}{2} \bigr)}, 
\end{equation*}
where in this case $\theta=\dfrac{\beta}{2(1-s)+\beta}\in (0,1)$. By \eqref{st3.40DFM}, we thus have

\begin{align} \nonumber
\Vert Du_j\Vert_{L^{\tilde{p}}(B_{\tau_1})}
&\leq 
\dfrac{c}{(\tau_2-\tau_1)^{1+\delta+\kappa}} [\mathbb{\tilde{F}}(u,B_r)]^\frac{1-\theta}{p}[\mathbb{F}(u,B_r)]^\frac{\theta}{2d}\\ \label{st-4p12}
&\qquad+ \dfrac{c}{(\tau_2-\tau_1)^
{1+\delta+\kappa}} 
[\mathbb{F}(u,B_r)]^\frac{\theta}{2d}
\Vert Du_j\Vert^{\left(p+\frac{\xi}{1-\gamma}\right)\frac{1-\theta}{p}}_{L^{p+\frac{\xi}{1-\gamma}}(B_{\tau_2})} \notag \\
& \qquad {+ c [\mathbb{F}(u,B_r)]^\frac{\theta}{2d}}.
\end{align}

\noindent By \eqref{controlpptil}, we are able to exploit the interpolation inequality \eqref{interp-5}, where $\theta_1$ is defined in the same way, with the difference that the expression of $\tilde{p}$ is different owing to \eqref{ptil}.

Then, inserting \eqref{interp-5} in \eqref{st-4p12} implies
\begin{align*} \nonumber
\Vert Du_j\Vert_{L^{\tilde{p}}(B_{\tau_1})}
&\leq 
c [\mathbb{\tilde{F}}(u,B_r)]^\frac{1-\theta}{p}[\mathbb{F}(u,B_r)]^\frac{\theta}{2d}\\ \label{st4-5p12}
&\qquad+c
[\mathbb{F}(u,B_r)]^\frac{\theta}{2d}
\Vert Du_j\Vert^{\theta_1\left(p+\frac{\xi}{1-\gamma}\right)\frac{1-\theta}{p}}_{L^{\tilde{p}}(B_{\tau_2})}
\Vert Du_j\Vert^{(1-\theta_1)\left(p+\frac{\xi}{1-\gamma}\right)\frac{1-\theta}{p}}_{L^{p}(B_{\tau_2})}\\
& \qquad {+ c [\mathbb{F}(u,B_r)]^\frac{\theta}{2d}}.
\end{align*}
Now we observe that, by \eqref{teta1} and \eqref{estlt1p12} we have
\begin{align*}
\theta_1\left(p+\frac{\xi}{1-\gamma}\right)\frac{1-\theta}{p} = 
\dfrac{\tilde{p}\xi}{\tilde{p}-p}\dfrac{1}{1-\gamma}\dfrac{1-s}{p(1-s)+\beta}<1.
\end{align*}
Applying Young's inequality and reabsorbing the term $\Vert Du_j\Vert_{L^{\tilde{p}}}$, we get
\begin{equation*}
\Vert Du_j\Vert_{L^{\tilde{p}}(B_{\tau_1})}\leq c[\mathbb{F}(u,B_r)]^{\kappa_1}[\mathbb{\tilde{F}}(u,B_r)]^{\kappa_2}.
\end{equation*}
Then, from \eqref{convtou} and weak lower semicontinuity, letting  $j\to +\infty$ we can conclude that
\begin{equation*}
Du \in L^{\tilde{q}}_{\textrm{loc}}(\Omega), \quad \text{for every} \ q<\tilde{q}<p+\dfrac{\min\lbrace \alpha,2\gamma\rbrace}{\Theta(1-\gamma)}.
\end{equation*}
This completes the proof of the theorem.

\subsection{Higher differentiability }
In this section, we give the proof of Theorem \ref{thmhd} taking advantage of the higher integrability result proved in Section \ref{sechi}.
According to Theorem \ref{thmhi}, we deduce that
\begin{equation}\label{gradxi}
Du \in L_{loc}^{p+\frac{\xi}{1-\gamma}}(\Omega, \mathbb{R}^{N\times n}),
\end{equation}
where $\xi$ was {introduced} in equation \eqref{sigma}, hence {$q<p+\frac{\xi}{1-\gamma}$.}
This implies that we do not need the approximating problems to study the higher fractional differentiability of $Du$.
\\

\noindent \textit{Proof of Theorem \ref{thmhd}.} {According to \eqref{gradxi}, a solution to problem \eqref{minprob} satisfies $u\in W^{1,q}_{loc}(\Omega, R^N)$, therefore $u$} solves the equation
\begin{equation}
\int_{B_R}\partial_zF(x,Du)\cdot D\fhi \ \mathrm{d}x=0, \quad \forall \varphi \in W^{1,q}_0(B_R,\R^N). \label{equation}
\end{equation}
Let us fix a ball $B_R$ such that $B_{2R} \Subset \Omega$ and a cut-off function $\eta \in \mathcal{C}_0^1(B_R)$ such that $0 \leq \eta \leq 1$, $\eta =1$ on $B_{R/2}$ and $|D \eta | \leq \frac{C}{R}$.
\\Now, for $|h| \leq R/4$, we consider test functions
\begin{equation}
\varphi_{1}(x)= \eta^{2}(x) \tau_hu(x)\label{test1}
\end{equation}
and
\begin{gather}
\varphi_{2}(x)= \eta^{2}(x-h) \tau_{-h}u(x). \label{test2}
\end{gather}
\\Inserting \eqref{test1} and \eqref{test2} in \eqref{equation}, we obtain
\begin{align}
\displaystyle\int_{\Omega}   \partial_zF(x,Du(x)) \cdot & D(\eta^2(x) \tau_hu(x))  \ \mathrm{d}x + \displaystyle\int_{\Omega}  \partial_zF(x,Du(x)) \cdot D(\eta^2(x-h) \tau_{-h}u(x)) \ \mathrm{d}x =0 .\label{2:4}
\end{align}
By means of a simple change of variable, we can write the second integral on the left hand side of the previous inequality as follows
\begin{align}
\displaystyle\int_{\Omega} \partial_zF(x+h,Du(x+h)) \cdot D(-\eta^2(x) \tau_hu(x))  \ \mathrm{d}x \label{2:5}
\end{align}
and so inequality \eqref{2:4} becomes
\begin{align}
\displaystyle\int_{\Omega} [\partial_zF(x+h,Du(x+h))-\partial_zF(x,Du(x))] \cdot D(\eta^2 \tau_hu(x))  \ \mathrm{d}x 
= 0 \label{2:6}
\end{align}
We can write previous inequality as follows
\begin{align}
0 = & \displaystyle\int_{\Omega} [\partial_zF(x+h,Du(x+h))-\partial_zF(x+h,Du(x))] \cdot \eta^{2}(x)D\tau_{h}u(x)  \ \mathrm{d}x\notag\\
  &+\displaystyle\int_{\Omega}[  \partial_zF(x+h,Du(x+h))-\partial_zF(x+h,Du(x)) ]\cdot 2\eta(x) D \eta(x)\tau_{h}u(x)) \ \mathrm{d}x\notag\\
  &+\displaystyle\int_{\Omega}[ \partial_zF(x+h,Du(x))-\partial_zF(x,Du(x))] \cdot \eta^{2}(x)D\tau_{h}u(x) \ \mathrm{d}x\notag\\
  &+\displaystyle\int_{\Omega}[\partial_zF(x+h,Du(x))-\partial_zF(x,Du(x))] \cdot 2\eta(x) D \eta(x)\tau_{h}u(x))  \ \mathrm{d}x \notag\\
 =:& I_{1}+I_{2}+I_{3}+I_{4}, \label{2:7}
\end{align}
that yields
\begin{align}
I_1 \leq & |I_2| + |I_3| + |I_4|  . \label{2:8}
\end{align}
The ellipticity assumption (F2) and Lemma \ref{D1} imply
\begin{align}
I_{1} \geq &\nu \displaystyle\int_{\Omega}  \eta^{2} (x)|\tau_{h}Du|^{2}(\mu^2 + |Du(x+h)|^{2}+|Du(x)|^{2})^{\frac{p-2}{2}} \ \mathrm{d}x \notag\\
\geq & C(\nu) \displaystyle\int_{B_{R/2}}   |\tau_{h}V_p(Du)|^{2} \ \mathrm{d}x ,\label{I1}
\end{align}
where in the last inequality we used that $\eta=1$ on the ball $B_{R/2}$.\\
From the growth condition (F4), H\"{o}lder's inequality and Lemma \ref{ldiff}, we get
\begin{align}
|I_{2}|  \leq & 2L \displaystyle\int_{\Omega} |D \eta(x)| \eta(x)  (1 + |Du(x+h)|^{2}+|Du(x)|^{2})^{\frac{q-1}{2}}|\tau_{h}u(x)| \ \mathrm{d}x \notag\\
\le
&\dfrac{C(L)}{R} \biggl( \displaystyle\int_{B_{R}}|\tau_h u(x)|^{q} \ \mathrm{d}x \biggr)^{\frac{1}{q}}  \biggl( \displaystyle\int_{B_{2R}}(1+|D u|)^{q} \ \mathrm{d}x \biggr)^{{\frac{q-1}{q}}} \notag \\
 \leq
&\dfrac{C(L)}{R} |h| \displaystyle\int_{B_{2R}}(1+|D u|)^q \ \mathrm{d}x .
\label{I3}
\end{align}
In order to estimate the integral $I_3$, we use assumption (F3), H\"{o}lder's inequality as follows
\begin{align}
|I_3| \leq & \displaystyle\int_{\Omega} \eta^2 |\tau_h Du| |h|^{\alpha} (g(x+h)+g(x))(1+|Du(x)|^2)^{\frac{q-1}{2}} \ \mathrm{d}x \notag\\
\le 
& |h|^{\alpha} \biggl( \displaystyle\int_{B_{2R}}g^\sigma \ \mathrm{d}x\biggr)^\frac{1}{\sigma} \biggl(  \displaystyle\int_{B_R}| \tau_h Du|^\frac{\sigma}{\sigma-1}(1+|Du|)^\frac{(q-1)\sigma}{\sigma-1} \ \mathrm{d}x \biggr)^{\frac{\sigma-1}{\sigma}}\notag\\
\le 
& |h|^{\alpha} \biggl( \displaystyle\int_{B_{2R}}g^\sigma \ \mathrm{d}x \biggr)^\frac{1}{\sigma} \biggl(  \displaystyle\int_{B_{2R}}(1+|Du|)^\frac{q\sigma}{\sigma-1} \ \mathrm{d}x\biggr)^{\frac{\sigma-1}{\sigma}}. \label{I4}
\end{align}
We note that
$$\dfrac{q\sigma}{\sigma-1}= p+\dfrac{\xi}{1-\gamma},$$
which ensures that the second integral on the right hand side of the previous inequality is finite.

Arguing analogously, we infer the following estimate for the integral $I_4$. 
\begin{align}
|I_4| \leq &\frac{C}{R} |h|^\alpha \int_{B_R}|\tau_h u| (g(x+h)+g(x))(1+|Du|^2)^\frac{q-1}{2} \ \mathrm{d}x \notag\\
\le 
& \frac{C}{R}|h|^{\alpha} \biggl( \displaystyle\int_{B_{2R}}g^\sigma \ \mathrm{d}x \biggr)^\frac{1}{\sigma} \biggl(  \displaystyle\int_{B_R}| \tau_h u|^\frac{q\sigma}{\sigma-1} \ \mathrm{d}x \biggr)^\frac{\sigma-1}{q\sigma}
 \biggl(  \displaystyle\int_{B_R}(1+|Du|)^\frac{q\sigma}{\sigma-1}dx\biggr)^{\frac{(q-1)(\sigma-1)}{q\sigma}}\notag\\
\le 
& \frac{C}{R}|h|^{\alpha+1} \biggl( \displaystyle\int_{B_{2R}}g^\sigma \ \mathrm{d}x \biggr)^\frac{1}{\sigma} \biggl(  \displaystyle\int_{B_{2R}}(1+|Du|)^\frac{q\sigma}{\sigma-1} \ \mathrm{d}x \biggr)^{\frac{\sigma-1}{\sigma}}
\label{I6}. 
\end{align}
Inserting estimates \eqref{I1},  \eqref{I3}, \eqref{I4} and \eqref{I6} in \eqref{2:8}, we infer

\begin{align}
     C & (\nu) \displaystyle\int_{B_{R/2}}   |\tau_{h}V_p(Du)|^{2} \ \mathrm{d}x \notag\\
    \le & \dfrac{C(L)}{R} |h| \displaystyle\int_{B_{2R}}(1+|D u|)^q \ \mathrm{d}x \notag\\
    &+|h|^{\alpha} \biggl( \displaystyle\int_{B_{2R}}g^\sigma \ \mathrm{d}x \biggr)^\frac{1}{\sigma} \biggl(  \displaystyle\int_{B_{2R}}(1+|Du|)^\frac{q\sigma}{\sigma-1} \ \mathrm{d}x \biggr)^{\frac{\sigma-1}{\sigma}} \notag\\
    &+ \frac{C}{R}|h|^{\alpha+1} \biggl( \displaystyle\int_{B_{2R}}g^\sigma \ \mathrm{d}x \biggr)^\frac{1}{\sigma} \biggl(  \displaystyle\int_{B_{2R}}(1+|Du|)^\frac{q\sigma}{\sigma-1} \ \mathrm{d}x \biggr)^{\frac{\sigma-1}{\sigma}},
\end{align}
that yields
\begin{align}
 \displaystyle\int_{B_{R/2}}  |\tau_{h}V_p(Du)|^{2} \ \mathrm{d}x 
    \le C |h|^\alpha \biggl\{
    \displaystyle\int_{B_{2R}}g^\sigma \ \mathrm{d}x + \displaystyle\int_{B_{2R}}
     (1+|Du|)^{p+\frac{\xi}{1-\gamma}} \ \mathrm{d}x
    \biggr\},
\end{align}
for a constant $C:=C(p,q,\gamma,\sigma,L,\nu,R)$.
\endproof

\section*{Acknowledgements}
The authors would like to thank Prof.\ Giuseppe Mingione for fruitful discussion. {The authors thank the anonymous referee for the careful reading of the manuscript and
for valuable comments.}
The authors are supported by GNAMPA (Gruppo Nazionale per l'Analisi Matematica, la Probabilit\`{a} e le loro Applicazioni) of INdAM (Istituto Nazionale di Alta Matematica).
 A.G. Grimaldi has been partially supported by the Gruppo Nazionale per l’Analisi Matematica, la Probabilità e le loro Applicazioni (GNAMPA) of the Istituto Nazionale di Alta Matematica (INdAM) through INdAM-GNAMPA
project (CUP\_E53C22001930001).


\begin{thebibliography}{}

\bibitem{1dfm}
E. Acerbi, G. Bouchitt\'{e} and I. Fonseca, Relaxation of convex functionals: the gap problem, \textit{Ann. Inst. H. Poincar\'{e} Anal. Non Lin\'{e}aire} \textbf{20} (2003), 359--390 


\bibitem{ambrosio1}
\newblock P. Ambrosio,
\newblock Besov regularity for a class of singular or degenerate elliptic equations,
\newblock \textit{J. Math. anal. Appl.} \textbf{505}(2) (2022)


\bibitem{ambrosio2}
\newblock P. Ambrosio,
\newblock Fractional Sobolev regularity for solutions to a strongly degenerate parabolic equation,
\newblock \textit{Forum Math.}. https://doi.org/10.1515/forum-2022-0293

\bibitem{baison.clop2017}
\newblock A. L. Baison, A. Clop, R. Giova, J. Orbitg and A. Passarelli di Napoli,
\newblock Fractional differentiability for solutions of nonlinear elliptic equations,
\newblock \textit{Potential Anal.} \textbf{46} (3) (2017), 403--430 

\bibitem{2dfm} A. Balci, L. Diening and M. Surnachev, New examples on Lavrentiev gap using fractals, \textit{Calc. Var. PDE}, \textbf{59} (2020) 180

\bibitem{BCM}
\newblock P. Baroni, M. Colombo and G. Mingione, \newblock Regularity for general functionals with double phase,
\newblock \textit{Calc.
Var. PDE,} \textbf{57} (2018) 62 


\bibitem{byun}
\newblock S.S. Byun, Y. Cho and J. Ok,
\newblock Global gradient estimates for nonlinear obstacle problems with non-standard growth,
\newblock \textit{Forum Math.} \textbf{28}(4) (2016), 729--747


\bibitem{byun1}
\newblock S.s. Byun and C. Namkyeong,
\newblock Higher differentiability for solutions of a general class of nonlinear elliptic obstacle problems with Orlicz growth,
\newblock \textit{NoDEA} \textbf{29}(73) (2022).  https://doi.org/10.1007/s00030-022-00807-x

\bibitem{carozza0}
\newblock M. Carozza, J. Kristensen, and A. Passarelli di Napoli,
\newblock Higher differentiability of minimizers of convex
variational integrals,
\newblock \textit{Ann. Inst. H. Poincaré Anal. Non Linéaire,} \textbf{28} (2011), 395--411 

\bibitem{carozza1}
\newblock M. Carozza, J. Kristensen, and A. Passarelli di Napoli,
\newblock Regularity of minimizers of autonomous convex variational integrals,
\newblock \textit{Ann. Sc. Norm. Super. Pisa Cl. Sci.} \textbf{13}(4) (2014), 1065--1089 


\bibitem{carozza2}
\newblock M. Carozza, J. Kristensen, and A. Passarelli di Napoli,
\newblock On the validity of the Euler-Lagrange system,
\newblock \textit{Comm. Pure Appl. Anal.}, \textbf{14}(1) (2018), 51--62 

\bibitem{caselli}
\newblock M. Caselli, A. Gentile and R. Giova,
\newblock Regularity results for solutions to obstacle problems with Sobolev coefficients,
\newblock \textit{J. Differential Equations} 269, 8308--8330 (2020)


\bibitem{defilippis1}
\newblock I. Chlebicka and C. De Filippis, \newblock Removable sets in non-uniformly elliptic problems,
\newblock \textit{Annali di Matematica
Pura ed Applicata} 199 (2) (2020), 619-649


\bibitem{choe}
H.J. Choe, Interior behaviour of minimizers for certain functionals with nonstandard growth, \textit{Nonlinear
Anal.} \textbf{19} (1992), 933--945 

\bibitem{clop}
\newblock A. Clop, R. Giova and A. Passarelli di Napoli,
\newblock Besov regularity for solutions of $p$-harmonic equations,
\newblock \textit{Adv. Nonlinear Anal.} (2017). https://doi.org/10.1515/anona-2017-0030

\bibitem{CM}
\newblock M. Colombo and  G. Mingione, 
\newblock Bounded minimisers of double phase variational integrals,
\newblock \textit{Arch. Ration. Mech. Anal.} 218, 219--273 (2015)


\bibitem{cupini.marcellini.mascolo.passarelli}
\newblock G. Cupini, P. Marcellini, E. Mascolo and A. Passarelli di Napoli,
\newblock Lipschitz regularity for degenerate elliptic integrals with $p,q$-growth,
\newblock \textit{Advances in Calculus of Variations} (2021). https://doi.org/10.1515/acv-2020-0120 


\bibitem{defilippis}
\newblock C. De Filippis,
\newblock Regularity results for a class of non-autonomous obstacle
problems with $(p,q)$-growth,
\newblock \textit{J. Math. Anal. Appl.} 501 (2021) 123450


\bibitem{20dfm}
C. De Filippis and G. Mingione, On the regularity of minima of non-autonomous functionals, \textit{J. Geom.
Anal.}, \textbf{30 }(2020), 1584--1626 


\bibitem{DFM}
\newblock C. De Filippis and G. Mingione, \newblock Interpolative gap bounds for nonautonomous integrals,
\newblock \textit{Analysis and Mathematical Physics}, \textbf{11} (2021), 117
https://doi.org/10.1007/s13324-021-00534-z.


\bibitem{eleuteri2007}
\newblock M. Eleuteri,
\newblock Regularity results for a class of obstacle problems,
\newblock \textit{Applications of Mathematics}, \textbf{52} (2) (2007), 137--170 


\bibitem{eleuteri.passarelli}
\newblock M. Eleuteri, A. Passarelli di Napoli, \newblock Higher differentiability for solutions to a class of obstacle problems,
\newblock \textit{Calc. Var.}, \textbf{57} (2018), 115  https://doi.org/10.1007/s00526-018-1387-x





\bibitem{eleuteri.passarelli1}
\newblock M. Eleuteri and A. Passarelli di Napoli,
\newblock Regularity results for a class of non-differentiable obstacle problems,
\newblock \textit{Nonlinear Analysis} 194 (2020) 111434




\bibitem{eleuteri.passarelli2021}
\newblock M. Eleuteri and A. Passarelli di Napoli,
\newblock On the validity of variational inequalities for obstacle problems with non-standard growth,
\newblock \textit{Annales Fennici Mathematici}, \textbf{47} (1) (2021), 395--416  https://doi.org/10.54330/afm.114655 

\bibitem{25dfm}
L. Esposito, F. Leonetti, G. Mingione, Sharp regularity for functionals with $(p,q)$ growth,\textit{ J. Differ.
Equ.}, \textbf{204} (2004), 5--55 



\bibitem{FMM}
\newblock I. Fonseca, J. Maly and G. Mingione, 
\newblock Scalar minimizers with fractal singular sets,
\newblock \textit{Arch. Ration. Mech.
Anal.} 172, 295–307 (2004)


\bibitem{gavioli1}
\newblock C. Gavioli,
\newblock A priori estimates for solutions to a class of obstacle problems under $p,q$-growth
conditions, 
\newblock \textit{Journal of Elliptic and Parabolic Equations} 5 (2) (2019), 325-347

\bibitem{gavioli2}
\newblock C. Gavioli,
\newblock Higher differentiability of solutions to a class of obstacle problems under non-standard
growth conditions,
\newblock \textit{Forum Mathematicum} (2019), 31(6), 1501--1516

\bibitem{Gentile}
\newblock A. Gentile,
\newblock Higher differentiability results for solutions to a class of non-autonomous obstacle problems with sub-quadratic growth conditions,
\newblock \textit{Forum Mathematicum},\textbf{ 33 }(3) (2021), 669--695 

\bibitem{gentile0}
\newblock A. Gentile and R. Giova,
\newblock Regularity results for solutions to a class of
non-autonomous obstacle problems with
sub-quadratic growth conditions,
\newblock Preprint arXiv:2201.07679v1




\bibitem{gentile1}
\newblock A. Gentile, R. Giova and A. Torricelli,
\newblock Regularity results for bounded solutions to obstacle
problems with non-standard growth conditions,
\newblock Preprint arXiv:2110.09586v1



\bibitem{Giova}
\newblock R. Giova,
\newblock Higher differentiability for $n$-harmonic systems with Sobolev coefficients,
\newblock \textit{J. Differ. Equ.} \textbf{259}(11) (2015), 5667--5687 

\bibitem{Giova.Pass}
\newblock R. Giova and A. Passarelli di Napoli, 
\newblock Regularity results for a priori bounded minimizers of non-autonomous functionals with discontinuous coefficients,
\newblock \textit{Adv. Calc. Var.} (2018). https://doi.org/10.1515/acv-2016

\bibitem{giusti}
\newblock E. Giusti, 
\newblock Direct methods in the calculus of variations,
\newblock \textit{World scientific publishing Co., Singapore} (2003) 


\bibitem{grimaldi0}
\newblock A.G. Grimaldi,
\newblock Regularity results for solutions to a class of obstacle problems,
\newblock \textit{Nonlinear Analysis: Real World Applications}, \textbf{62} (2021), 103377

\bibitem{grimaldi}
\newblock A.G. Grimaldi,
\newblock \textit{Higher differentiability for bounded solutions to a class of obstacle problems with $(p,q)$-growth},
\newblock \textit{Forum Mathematicum}, \textbf{35} (2) (2023), 457--485 


\bibitem{grimaldi.ipocoana}
\newblock A.G. Grimaldi and E. Ipocoana,
\newblock Higher fractional differentiability for solutions to a class of obstacle problems with non-standard growth conditions,
\newblock \textit{Advances in Calculus of Variations}, DOI: 10.1515/acv-2021-0074

\bibitem{grimaldi.ipocoana1}
\newblock A.G. Grimaldi and E. Ipocoana,
\newblock Higher differentiability results in the scale of Besov spaces to a class of double-phase obstacle problems,
\newblock \textit{ESAIM: COCV} 28 (2022) 51, doi: 10.1051/cocv/2022050


\bibitem{haroske} D. Haroske, 
\newblock Envelopes and sharp embeddings of function spaces,
\newblock \textit{Chapman and Hall CRC, Boca Raton} (2006)


\bibitem{koskela}
\newblock P. Koskela, D. Yang and Y. Zhou,
\newblock Pointwise characterizations of Besov and Triebel-Lizorkin spaces and quasiconformal mappings,
\newblock \textit{Adv. Math},\textbf{ 226} (4) (2011), 3579--3621 


\bibitem{marcellini1986} 
\newblock P. Marcellini, 
\newblock On the definition and the lower semicontinuity of certain quasiconvex integrals,
\newblock \textit{Ann. Inst. H. Poincaré Anal. Non Linéaire}, \textbf{3}(5) (1986), 391--409

\bibitem{marcellini1989} 
P. Marcellini, Regularity of minimizers of integrals of the calculus of variations with nonstandard growth conditions,
\textit{Arch. Ration. Mech. Anal.} \textbf{105}(3) (1989), 267--284

\bibitem{marcellini1991} 
\newblock P. Marcellini, 
\newblock Regularity and existence of solutions of elliptic equations with $p,q$-growth conditions,
\newblock \textit{J. Differential Equations}, \textbf{90}(1) (1991), 1-30



\bibitem{Pass1}
\newblock A. Passarelli di Napoli,
\newblock Higher differentiability of minimizers of variational integrals with Sobolev coefficients,
\newblock \textit{Adv. Calc. Var.} \textbf{7}(1) (2014), 59--89

\bibitem{Pass2}
\newblock A. Passarelli di Napoli, 
\newblock Higher differentiability of solutions of elliptic systems with Sobolev coefficients: the case $p=n=2$,
\newblock \textit{Potential Anal.} \textbf{41}(3) (2014), 715--735 

\bibitem{Pass3}
\newblock A. Passarelli di Napoli,
\newblock Regularity results for non-autonomous variational integrals with discontinuous coefficients,
\newblock \textit{Atti Accad. Naz. Lincei Rend. Lincei Mat. Appl} \textbf{26}(4) (2015), 475--496 






\bibitem{triebel83}
\newblock H. Triebel,
\newblock Theory of function spaces,
\newblock \textit{Birkh\"{a}user-Verlag}, Basel (1983)

\bibitem{zhang.zheng}
\newblock X. Zhang and S. Zheng,
\newblock Besov regularity for the gradients of solutions to non-uniformly elliptic obstacle problems,
\newblock \textit{Journal of Mathematical Analysis and Applications} 505 (2) (2021)

\bibitem{60dfm}
\newblock V.V. Zhikov, 
\newblock On Lavrentiev’s phenomenon,
\newblock \textit{ Russ. J. Math. Phys. }\textbf{3}, 249–269 (1995)

\bibitem{61dfm}
\newblock V.V. Zhikov,
\newblock On some variational problems,
\newblock \textit{Russ. J. Math. Phys.} \textbf{5}, 105–116 (1997)

\bibitem{62dfm}
\newblock V.V. Zhikov, 
\newblock Lavrentiev phenomenon and homogenization for some variational problems,
\newblock \textit{C. R. Acad.
Sci. Paris S\'{e}r. I Math.} \textbf{316}, 435–439 (1993)
\end{thebibliography}
\end{document}